\input amstex.tex
\documentstyle{amsppt}
\input pictex.tex
\magnification=\magstep1
\pagewidth{15.9truecm}
\nologo
\NoRunningHeads
\def\Pic{\operatorname{Pic}}

\def\supp{\operatorname{supp}}
\def\Spec{\operatorname{Spec}}
\def\sing{\operatorname{sing}}
\def\ord{\operatorname{ord}}
\def\div{\operatorname{div}}
\def\dim{\operatorname{dim}}

\def\card{\operatorname{card}}
\refstyle{A}
\widestnumber\key{AKMW}
\topmatter
\title
Zeta functions and \lq Kontsevich invariants\rq\ on singular varieties
\endtitle
\author
Willem Veys  \bigskip
\endauthor
\address K.U.Leuven, Departement Wiskunde, Celestijnenlaan 200B,
         B--3001 Leuven, Belgium  \endaddress
\email wim.veys\@wis.kuleuven.ac.be  \newline 
{}\quad
http://www.wis.kuleuven.ac.be/wis/algebra/veys.htm
\endemail
\abstract
Let $X$ be a nonsingular algebraic variety in characteristic zero.
To an effective divisor on $X$ Kontsevich has associated a certain motivic integral, living in a completion of the Grotendieck ring of algebraic varieties. He used this invariant to show that birational (smooth, projective) Calabi--Yau varieties have the same Hodge numbers.
Then Denef and Loeser introduced the invariant {\it motivic (Igusa) zeta function}, associated to a regular function on $X$, which specializes to both the classical $p$--adic Igusa zeta function and the topological zeta function, and also to Kontsevich's invariant.

This paper treats a generalization to singular varieties. Batyrev already considered such a \lq Kontsevich invariant\rq\ for log terminal varieties (on the level of Hodge polynomials of varieties instead of in the Grothendieck ring), and previously we introduced a motivic zeta function on normal surface germs. Here on any $\Bbb Q$--Gorenstein variety $X$ we associate a motivic zeta function and a 
\lq Kontsevich invariant\rq\ to effective $\Bbb Q$--Cartier divisors on $X$ whose support contains the singular locus of $X$.
\endabstract
\endtopmatter
\bigskip
\document
\head
Introduction
\endhead
\bigskip
\noindent
{\bf 0.1.} Let $k$ be a field of characteristic zero.  To a nonsingular
(irreducible) variety $X$ and a morphism $f : X \rightarrow \Bbb A^1$, both defined over $k$, was
associated the invariant {\it motivic (Igusa) zeta function} by Denef and
Loeser
[DL2].  By definition it lives in a power series ring in one variable over the
ring $\Cal M_L$, where $\Cal M$ is the Grothendieck ring of algebraic varieties
over $k$, $L$ is the class of $\Bbb A^1$ in $\Cal M$, and $\Cal M_L$ denotes
localization.  When $X = \Bbb A^d$ this invariant specializes to both the usual
$p$--adic Igusa zeta function and the topological zeta function associated to a
polynomial $f$.  (In fact in [DL2] the authors treat an even more general
invariant, involving motives instead of varieties, from which also the whole
Hodge spectrum of $f$ at any point of $f^{-1} \{ 0 \}$ can be deduced.)  This
notion of motivic zeta function can easily be extended to an effective divisor
$D$ instead of just a morphism $f$.

The authors were inspired by Kontsevich's idea of {\it motivic integration}.
In [Kon] Kontsevich associated to a nonsingular irreducible variety $X$ and an
effective divisor $D$ on $X$ an invariant $\Cal E (D)$, living by definition
in an appropriate completion $\hat \Cal M$ of $\Cal M_L$.  He used this
invariant to show that birationally equivalent (smooth, projective) Calabi--Yau
varieties have the same Hodge numbers.
\bigskip
\noindent
{\bf 0.2.} There are important formulas for these invariants in terms of an
embedded resolution (with strict normal crossings) $h : Y \rightarrow X$ of
supp $D$.  Let $\dim X = d$ and denote by $E_i, i \in T$, the irreducible
components of $h^{-1} (\supp D)$.  To the $E_i$ are associated natural
multiplicities $N_i$ and $\nu_i$ defined by $h^\ast D = \sum_{i \in T} N_i E_i$
and $\div(h^\ast dx) = \sum_{i \in T} (\nu_i - 1)E_i$, where $dx$ is a local
generator of the sheaf of regular differential $d$-forms on $X$.  Also we
partition $Y$ into the locally closed strata 
$E^\circ_I := (\cap_{i \in I} E_i) \setminus (\cup_{\ell \not\in I} E_\ell)
\text
{ for } I \subset T.$ 

We denote the class of a variety $V$ in $\Cal M$ by $[V]$, and by analogy with
the usual $p$--adic Igusa zeta function we denote the variable of the power
series ring over $\Cal M_L$ formally by $L^{-s}$.  Then the motivic zeta
function $\Cal Z(D,s)$ of $D$ is given by the formula
$$\Cal Z(D,s) = L^{-d} \sum_{I \subset T} [E^\circ_I] \prod_{i \in I}
\frac{(L-1)L^{-\nu_i}(L^{-s})^{N_i}}{1 - L^{-\nu_i} (L^{-s})^{N_i}}$$
and so it lives already in a localization of the polynomial ring $\Cal
M_L[L^{-s}]$.  Kontsevich's invariant for $D$ is given by 
$$\Cal E(D) = L^{-d} \sum_{I \subset T} [E^\circ_I] \prod_{i \in I} \frac{L -
1}{L^{\nu_i+N_i} - 1}$$
and can thus in some sense be derived from $\Cal Z(D,s)$ by `substituting
$s=1$'.
\bigskip
\noindent
{\bf 0.3.} One can specialize $\Cal Z(D,s)$ and $\Cal E(D)$ to more
`concrete' invariants, involving instead of the class $[V]$ of a variety $V$ in
$\Cal M$ other additive invariants as the Hodge polynomial $H(V)$ or the Euler
characteristic $\chi(V)$ of $V$.  With a little work one obtains  for instance from $\Cal Z (D,s)$ the {\it
topological zeta function}
$$z(D,s) = \sum_{I \subset T} \chi (E^\circ_I) \prod_{i \in I} \frac{1}{\nu_i +
sN_i} \in \Bbb Q(s) \, ,$$
which was introduced in [DL1] for $X=\Bbb A^d$ and $k=\Bbb C$,
and the invariant
$$e(D) = \sum_{I \subset T} \chi (E^\circ_I) \prod_{i \in I} \frac{1}{\nu_i +
N_i} \in \Bbb Q \, .$$
\bigskip
\noindent
{\bf 0.4.}  Can the invariants above be generalized to singular (normal)
varieties $X$ such that analogous formulas in terms of an embedded resolution
are valid ?  The main problem is whether these formulas are independent of the
chosen resolution.  Let $D$ be an effective Weil divisor on $X$ and $h : Y
\rightarrow X$ an embedded resolution of $X_{\text{sing}} \cup \supp D$ with
irreducible components $E_i, i \in T$, of $h^{-1}(X_{\text{sing}} \cup \supp
D)$.  Can we generalize the multiplicities $N_i$ and $\nu_i$ ?  When $D$ is
Cartier (or $\Bbb Q$--Cartier) the same expression $h^\ast D = \sum_{i \in T}
N_i E_i$ makes sense.  We think that the most natural generalization of the
$\nu_i$ are the log discrepancies given by $K_Y = h^\ast K_X + \sum_{i \in T}
(\nu_i - 1)E_i$, where $K_{\cdot}$ is the canonical divisor. 
To this end we need in
general $X$ to be Gorenstein (or $\Bbb Q$--Gorenstein).  Up to now the following
generalizations appeared (with $k = \Bbb C$).

(a) In dimension 2 these multiplicities are defined for arbitrary Weil
divisors on normal surfaces.  In [V3] we introduced a topological zeta function
and a motivic zeta function for effective divisors on normal surface germs.  We
could have done this as well globally, associating to an effective Weil divisor
$D$ on a normal surface $X$ for which $X_{\text{sing}} \subset \supp D$
the zeta functions $\Cal Z(D,s)$ and $z(D,s)$, given by the same formulas as
above.

(b) In arbitrary dimension Batyrev [B2] considered the case $D = 0$ and
associated `Kontsevich--like' invariants to a log terminal $X$ on the level of
Hodge polynomials and Euler characteristics.  The last one, which he called
{\it stringy Euler number}, is given by the formula for $e(D)$ in (0.3) with
all $N_i =
0$.  The invariant on Hodge polynomial level was used in [B2] to define {\it
stringy Hodge numbers} for projective canonical Gorenstein varieties, and to
formulate a topological mirror duality test for canonical Calabi--Yau varieties.

(c) Batyrev [B3] also extended his construction to Kawamata log terminal pairs
$(X,D)$, i.e. pairs such that $K_X + D$ is $\Bbb Q$--Cartier and all $a_i > 0$
in the expression $K_Y = h^\ast (K_X + D) + \sum_{i \in T} (a_i - 1)E_i$.  On
the Euler characteristic level this invariant is given by the formula 
$$e \big((X,D)\big) = \sum_{I \subset T}  \chi (E^\circ_I) \prod_{i \in I}
\frac{1}{a_i}.$$
In [B3] these invariants are used to prove a version of Reid's McKay
correspondence conjecture.

We should mention that Batyrev is naturally restricted to the log terminality
conditions above (all $\nu_i > 0$ and all $a_i > 0$, respectively) by applying
motivic integration techniques to show that the formulas above are independent
of the chosen resolution; see [B2, Theorem 6.28].

\medskip
We also want to remark that $\Cal E(D)$ is generalized in [DL3] in a different way (see 3.5). 
\bigskip
\noindent
{\bf 0.5.} In this paper we extend the invariants above beyond the log
terminal case to the following general situation.  Now let $X$ be any normal
$\Bbb Q$--Gorenstein variety and $D$ an effective $\Bbb Q$--Cartier divisor with
$X_{\text{sing}} \subset \supp D$.
We associate first to these data zeta functions $\Cal Z(D,s), Z(D,s)$ and
$z(D,s)$ on `motivic' level, Hodge polynomial level and Euler characteristic
level, respectively, such that the same formulas as in (0.2) and (0.3) are
valid.  Then we {\it define} `Kontsevich' invariants $\Cal E(D), E(D)$ and
$e(D)$ on the analogous levels by taking the limit for $s \rightarrow 1$ in the
associated zeta functions (admitting the value $\infty$).  In particular when
all $\nu_i + N_i \ne 0$ the formulas in (0.2) and (0.3) are again valid.

Furthermore taking the limit for $s \rightarrow -1$ in the zeta functions we
obtain invariants $\Cal E \bigl((X,D) \bigr), E \bigl((X,D)\bigr)$ and $e \bigl((X,D)\bigr)$ of the pair $(X,D)$
on the same levels,
the last one given by the same formula as in (0.4).  

In fact we can relax our condition $X_{\text{sing}} \subset \supp D$ to
$LCS(X) \subset \supp D$, where $LCS(X)$ is the locus of log canonical
singularities of $X$.  In particular this locus is empty when $X$ is log
terminal; so we really generalize the invariants of [B2].
\bigskip
\noindent
{\bf 0.6.} In \S 1 we recall the motivic zeta function of Denef and
Loeser and the invariant of Kontsevich on smooth varieties $X$, generalizing
the first one to effective divisors instead of regular functions.  As an
introduction to singular varieties we treat the easy case of a 
canonical $X$ in \S 2; there we also consider an application to minimal models.
For $\Bbb Q$--Gorenstein varieties 
$X$ the zeta functions $Z(D,s)$ and $z(D,s)$ on the level of Hodge polynomials
and Euler characteristics, respectively, are constructed in an elementary way
in \S 3. We provide some examples in \S 4.  The \lq motivic\rq\ version requires more work.  In \S 5 we first introduce
a motivic zeta function $\Cal Z(D,J,s)$ on a smooth $X$, associated to both an
effective divisor $D$ and an invertible subsheaf $J$ of the sheaf of regular
differential forms on $X$.  (This can be compared with associating a $p$--adic
Igusa zeta function to both a polynomial and a differential form.)  Then we use
this object to define the motivic zeta function $\Cal Z(D,s)$ for  a $\Bbb
Q$--Gorenstein $X$ in \S 6.  We include an appendix indicating how to extend the
original Kontsevich invariant on smooth $X$ to $\Bbb Q$--divisors instead of
(ordinary) divisors, needing a finite extension of $\hat \Cal M$.

\bigskip
\noindent
{\bf 0.7.} {\sl Remark.}  After this work was finished we learned about the
proofs of
W{\l}odarczyk [W{\l}] and of Abramovich et al [AKMW] of the weak factorization
conjecture for birational maps.  Using weak factorization we can give another
proof that the zeta functions in this paper are well defined.
\bigskip
\bigskip
\head
1. Smooth varieties
\endhead
\bigskip
\noindent
{\bf 1.1.} Let $k$ be a field of characteristic zero; the varieties and
morphisms we
will consider are assumed to be defined over $k$.  (A variety is a reduced
separated scheme of finite type over $k$, not necessarily irreducible.)  

We fix some terminology concerning resolution.  A {\it resolution} of an
irreducible variety $X$ is a proper birational morphism $h : Y \rightarrow X$
from a smooth variety $Y$, which is an isomorphism outside the set $X_{\text
{sing}}$ of singular points of $X$.  A {\it log resolution} or {\it embedded
resolution} of an irreducible variety $X$ is a resolution $h : Y \rightarrow X$
of $X$ for which $h^{-1} (X_{\text {sing}})$ is a divisor with strict normal
crossings, i.e. with smooth irreducible components intersecting transversely.
A {\it log resolution} or {\it embedded resolution} of a reduced Weil divisor
$D$ on a normal variety $X$ is a proper birational morphism $h : Y \rightarrow
X$ from a smooth $Y$, which is an isomorphism outside $X_{\text{sing}} \cup D$,
and such that $h^{-1}(X_{\text{sing}} \cup D)$ is a divisor with strict
normal crossings.

We
denote by $\Cal M$ the Grothendieck ring of (algebraic) varieties
over $k$.  This is the free abelian group generated by the symbols $[V]$, where $[V]$ is a variety, subject to the relations  $[V] = [V^\prime]$ if $V \cong
V^\prime$ and $[V] = [V \setminus V^\prime] + [V^\prime]$  if $V^\prime$ is closed in $V$. Its ring structure is given by 
$[V] \cdot [V^\prime] := [V \times V^\prime]$.  We abbreviate $L := [\Bbb A^1]$
and denote by $\Cal M_L = \Cal M[L^{-1}]$ the localization of $\Cal M$ w.r.t.
the multiplicative set $\{ L^n, n \in \Bbb N \}$.  
\bigskip
\noindent
{\bf 1.2.}  For $[V] \in \Cal M$ we denote by $H(V) \in \Bbb Z[u,v]$ its Hodge
polynomial and by $\chi(V)$ its Euler characteristic.  We briefly explain
these notions.

\noindent
Let first $k = \Bbb C$.  Then for a variety $V$ we denote by
$h^{p,q} (H^i_c (V, \Bbb C))$ the rank of the $(p,q)$--Hodge component of its
$i$-th cohomology group with compact support and by $e^{p,q}(V) := \sum_{i \geq
0} (-1)^i
h^{p,q} (H^i_c(V,\Bbb C))$ its Hodge numbers.  The {\it Hodge
polynomial}  of $V$ is $H(V) = H(V;u,v) := \sum_{p,q} e^{p,q} (V) u^p v^q \in
\Bbb Z[u,v]$.

\noindent
Precisely by the defining relations of $\Cal M$ there is a well defined
ring morphism $H : \Cal M \rightarrow \Bbb Z[u,v]$ determined by $[V] \mapsto
H(V)$.

We denote by $\chi(V)$ the topological Euler characteristic of $V$, i.e.
the alternating sum of the ranks of its Betti or de Rham cohomology groups.
Clearly $\chi (V) = H(V;1,1)$ and we also obtain a ring morphism $\chi : \Cal M
\rightarrow \Bbb Z$ determined by $[V] \rightarrow \chi(V)$.

For arbitrary $k$ (of characteristic zero) we choose an embedding of the field of definition of the variety $V$ into $\Bbb C$.  Then we can define the same morphisms $H$ and
$\chi$ on $\Cal M$ starting from the $e^{p,q}(V)$; they are independent of the
chosen embedding since for a smooth projective $V$ we have that $e^{p,q}(V) =
(-1)^{p+q} \dim_k H^q(V,\Omega^p_V)$.
\bigskip
\noindent
{\bf 1.3.}  Till the end of this section we let $X$ be a smooth irreducible
variety
of dimension $d$ and $W$ a subvariety of $X$.

In [DL2] Denef and Loeser associate to $W \subset X$ and a morphism $f : X
\rightarrow \Bbb A^1$ an invariant named {\it motivic Igusa zeta function}.  We
recall here briefly its definition but generalize immediately to effective
divisors $D$ (instead of functions $f$).  We refer to [DL2] for more details
and motivation, and for the relation with the usual $p$--adic Igusa zeta
function.

We denote by $\Cal L(X)$ the {\it scheme of germs of arcs} on $X$.  It is a
scheme over $k$ whose $k$--rational points are the morphisms $\Spec  k[[t]]
\rightarrow X$ (called the germs of arcs on $X$).  In fact $\Cal L(X)$ is
defined as the projective limit $\varprojlim \Cal L_n(X)$ of the schemes of
truncated
arcs $\Cal L_n(X)$, whose $k$--rational points are the morphisms $\Spec (k[t]/t^{n+1}
k[t]) \rightarrow X$ (see [DL2] and [BLR, p.276]).  There are canonical
morphisms $\pi_n : \Cal L(X) \rightarrow \Cal L_n(X)$, induced by truncation.
Remark also that $\Cal L_0(X) = X$.

Now let $D$ be an effective divisor on $X$.  For $n \in \Bbb N$ we define
$Y_{n,D,W}$ as the subscheme of $\Cal L(X)$ whose $K$--rational points, for any
field $K \supset k$, are the morphisms $\varphi : \Spec K[[t]]
\rightarrow X$ satisfying the following conditions :

\itemitem{(i)} $\varphi$ sends the closed point of $\Spec K[[t]]$ to a point $P$
in $W$;
 	
\itemitem{(ii)} if $f$ is a local equation of $D$ at $P$, then the power series in $t$ given by $f \circ \varphi$ must be exactly of order $n$. 
(This is clearly independent of the choice of
$f$.)
\medskip
We then denote by $X_{n,D,W}$ the image of $Y_{n,D,W}$ in  $\Cal L_n (X)$,
viewed as a reduced subscheme of $\Cal L_n(X)$. 
The {\it motivic zeta function} of $D$ (and $W \subset X$) is
$$\Cal Z_W(D,s) = \Cal Z_W (X,D,s) := \sum_{n \in \Bbb N} [X_{n,D,W}]
L^{-(n+1)d-ns} \in \Cal M_L [[L^{-s}]].$$
Here $L^{-s}$ is just a variable and in the power series ring $\Cal M_L
[[L^{-s}]]$ we abbreviate $L^a \cdot (L^{-s})^b$ by $L^{a-sb}$ for  $a \in \Bbb
Z$ and $b \in \Bbb N$. 
(When $D$ is given by a global function $f$ on $X$ Denef and Loeser denoted
this invariant by $\int_W^\sim f^s$ in [DL2].)

One can think here mainly about $W$ as being $X$ itself, the divisor  
$\{f=0\}$, or a point of $\{f=0\}$. This $W$--formalism enables us to treat these cases together, and the greater generality is also useful.
\bigskip
\noindent
{\bf 1.3.1.}
We briefly compare this with the classical $p$--adic situation.  Let $f \in \Bbb Q_p[x] = \Bbb Q_p[x_1, \dots , x_d]$ and denote by $|z| = p^{-\ord_p z}$ the $p$--adic absolute value of $z \in \Bbb Q_p$.  {\it Igusa's local zeta function}
of $f$ is
$$Z_p(f,s) := \int_{\Bbb Z^d_p} |f(x)|^s |dx|$$
for $s \in \Bbb C$ with  $\Re(s) > 0$, where $|dx|$ denotes the Haar measure on $\Bbb
Q^d_p$ such that $\Bbb Z_p^d$ has measure $1$. When $f \in \Bbb Z_p[x]$ it is not difficult to verify that
$$Z_p(f,s) = \sum_{n \in \Bbb N} \card(X_{n,f}) p^{-(n+1)d-ns},$$
where $X_{n,f}$ is the image in $(\Bbb Z_p/p^{n+1} \Bbb Z_p)^d$ of $Y_{n,f} =
\{ x \in \Bbb Z^d_p | \ord_p f(x) = n \}$. See [D2] for an introduction and an overview on Igusa's local zeta function.
\bigskip
\noindent
{\bf 1.4.}  There is an important formula for $\Cal Z_W (D,s)$ in terms of a
log
resolution of $\supp D$.  In particular it implies the rationality
result that $\Cal Z_W (D,s)$ belongs in fact already to a certain localization
of the polynomial ring $\Cal M_L [L^{-s}]$.

Let $h : Y \rightarrow X$ be a log resolution of $\supp D$.  We denote by $E_i,
i \in T$, the irreducible components of $h^{-1}(\supp D)$ and by $N_i$ and
$\nu_i  - 1$ the multiplicities of $E_i$ in $h^\ast D$ and the divisor of
$h^\ast \text{dx}$, respectively, where $\text{dx}$
is a local generator of the sheaf $\Omega^d_X$ of regular differential
$d$--forms.  We
partition $Y$ into the locally closed strata $E^\circ_I := (\cap_{i \in I} E_i)
\setminus (\cup_{\ell \not\in I} E_\ell)$ for $I \subset T$.  (Here $E_\phi = Y
\setminus \cup_{\ell \in T} E_\ell$.)
\bigskip
\proclaim
{Theorem}  We have the formula
$$
\Cal Z_W(D,s) = L^{-d} \sum_{I \subset T}
[E^\circ_I
\cap
h^{-1}W] \prod_{i \in I} \frac{L-1}{L^{\nu_i + sN_i} - 1}
$$ 

\noindent
\text { (where}
  $\frac{1}{L^{\nu_i + sN_i}-1} := \frac{L^{-\nu_i-sN_i}}{1 - L^{-\nu_i -
sN_i}})$.  So $\Cal Z_W(D,s)$ belongs already to the localization 
$\Cal
M_L[L^{-s}]_{{{{(1 - L^{-n-Ns})}}_{{n,N \in \Bbb N \setminus \{ 0 \} }}}}$
of the polynomial ring $\Cal M_L [L^{-s}]$.
\endproclaim
\bigskip
\noindent
{\bf 1.4.1.} One should compare this formula with the classical formula of Denef [D1, Theorems 2.4 and 3.1] for Igusa's local zeta function $Z_p(f,s)$ of $f \in \Bbb Q[x_1,\dots,x_d]$ in terms of a resolution $h:Y \rightarrow \Bbb A^d$ of $\{f=0\}$. Using the notation above we have for all but finitely many $p$ that 
$$
Z_p(f,s) = p^{-d} \sum_{I \subset T}
\#(E^\circ_I)_{\Bbb F_p} \prod_{i \in I} \frac{p-1}{p^{\nu_i + sN_i} - 1} \, ,
$$ 
where $\#(\cdot)_{\Bbb F_p}$ denotes the number of $\Bbb F_p$--rational points of the reduction$\mod p$. See [D1, D2] for more details.

\bigskip
\noindent
{\bf 1.5.} Here we generalize $\Cal Z_W(D,s)$ to effective $\Bbb Q$--divisors on
$X$.
Now let $D$ be an effective $\Bbb Q$--divisor on $X$ and say that $rD$ is a
divisor for $r \in \Bbb N \setminus \{ 0 \}$.  We define $\Cal Z_W (D,s) :=
\Cal Z_W(rD,s/r)$, meaning by this the motivic zeta function of 1.3 for the
divisor $rD$, where the variable $L^{-s}$ is replaced by a variable
$(L^{-s})^{1/r}$.  This definition is easily checked to be independent of the
chosen $r$, using Theorem 1.4.\hfill\break
Moreover Theorem 1.4 is still valid in this context.  The only difference is
that the $N_i, i \in T$, are now rational numbers (of the form $a/r$ with $a
\in \Bbb N \setminus \{ 0 \})$, and one should consider $L^{-sN_i}$ as an
abbreviation of $((L^{-s})^{1/r})^{rN_i}$.
\bigskip
\noindent
{\bf 1.6.}  One can specialize the motivic zeta functions $\Cal Z_W (D,s)$
to more `concrete' invariants on the level of Hodge polynomials and on the
level of Euler characteristics.

($i$)
Let $D$ be an effective divisor on $X$.  Since the Hodge polynomial $H(\Bbb
A^1) = uv$ the morphism $H : \Cal M \rightarrow \Bbb Z[u,v]$ extends naturally
to a ring morphism $H : \Cal M_L \rightarrow \Bbb Z[u,v]_{uv} = \Bbb
Z[u,v][(uv)^{-1}]$ (and further to a morphism on power series rings over these
rings).  We define
$$Z_W(D,s) = Z_W(X,D,s) := H(\Cal Z_W(D,s)) = \sum_{n \in \Bbb N}
H(X_{n,D,W})(uv)^{-(n+1)d-ns},$$
where now we denote the variable of the power series ring over $\Bbb
Z[u,v]_{uv}$ by $(uv)^{-s}$.  Using the notation of 1.4 we have the formula
$$\aligned
Z_W(D,s) & = (uv)^{-d} \sum_{I \subset T} H(E^\circ_I \cap h^{-1}W) \prod_{i
\in
I} \frac{uv - 1}{(uv)^{\nu_i+sN_i}-1} \\
& \in \Bbb Z[u,v]_{uv} [(uv)^{-s}]_{(1-(uv)^{-n-Ns})_{n,N \in \Bbb N \setminus
\{ 0 \}}} \subset \Bbb Q(u,v)\big((uv)^{-s}\big).
\endaligned
$$

($ii$) To specialize further to the level of Euler characteristics one
takes heuristically the limit of the expression above for $u,v \rightarrow 1$.
We briefly explain the exact argument; see [DL2, (2.3)] for the argument
starting from $\Cal Z_W(D,s)$.  Let $R$ denote the subring of $\Bbb
Z[u,v]_{uv}[[(uv)^{-s}]]$ generated by $\Bbb Z[u,v]_{uv}[(uv)^{-s}]$ and the
elements $\frac{uv-1}{1 - (uv)^{-n-Ns}}$, where $n,N \in \Bbb N \setminus \{ 0
\}$.  $(Z_W(D,s)$ lives in $R$.)   By expanding $(uv)^{-s}$ and
$\frac{uv-1}{1-(uv)^{-n-Ns}}$ formally into series in $uv-1$, one constructs a
canonical algebra morphism 
$$R \rightarrow \Bbb Z[u,v]_{uv}[s][(n+sN)^{-1}]_{n,N \in \Bbb N \setminus \{ 0
\}} [[uv-1]],$$
where $[[uv-1]]$ denotes completion with respect to the ideal $(uv-1)$.
Composing this morphism with the quotient map given by dividing out $(uv-1)$ in
this last algebra yields a morphism
$$\varphi : R \rightarrow \frac{\Bbb Z[u,v]_{uv}}{(uv-1)} [s]
[(n+sN)^{-1}]_{n,N \in \Bbb N \setminus \{ 0 \}}.$$
In this last ring the evaluation $u=v=1$ is well defined; we put
$$\aligned
z_W(D,s) = z_W(X,D,s) & := \lim_{u,v \rightarrow 1} \varphi (Z_W(D,s)) \\
& = \sum_{I \subset T} \chi (E^\circ_I \cap h^{-1}W) \prod_{i \in I}
\frac{1}{\nu_i + sN_i} \in \Bbb Q(s).
\endaligned
$$
When $X = \Bbb A^n$ and $D$ is given by a polynomial $f$ these invariants are
just the {\it topological zeta functions} $Z_{\text{top}}(f,s)$ and
$Z_{\text{top},0}(f,s)$ of [DL1] if we take $W=X$ and $W = \{ 0 \}$,
respectively.

($iii$) As in 1.5 we can consider $Z_W (D,s)$ and $z_W(D,s)$ also for
$\Bbb Q$--divisors $D$.
\bigskip
\noindent
{\bf 1.7.}  Now we recall the original motivic integral, introduced by
Kontsevich
in [Kon], using the notation of 1.3.  We refer to [DL3] for a detailed
exposition in a much more general setting; see also the appendix. 
A nice introduction is [C].

We say that dim $M \leq n$ for $M \in \Cal M$ if $M$ can be expressed as a
$\Bbb Z$--linear combination of classes of algebraic varieties of dimension at
most $n$.  We consider the decreasing filtration $(F^m)_{m \in \Bbb Z}$ on
$\Cal M_L$, where $F^m$ is the subgroup of $\Cal M_L$ generated by $\{ [V]
L^{-i} | \dim V - i \leq -m \}$, and we denote by $\hat \Cal M$ the completion
of
$\Cal M_L$ with respect to this filtration.

Let again $D$ be an effective divisor on $X$.  We set 
$$\Cal E_W(D) = \Cal E_W (X,D) := \sum_{n \in \Bbb N}
\frac{[X_{n,D,W}]}{L^{(n+1)d}} L^{-n} \in \hat \Cal M ;$$
this expression converges in $\hat \Cal M$ since $\dim [X_{n,D,W}] \leq (n+1)d$.
This invariant was denoted as $[ \int_X e^D]$ by Kontsevich (for $W = X$) and
as
$\int_{\pi^{-1}_0 W} L^{-\ord_t \Cal O(-D)} d \mu$ in [DL3].  In this last
paper
Denef and Loeser develop an integration theory for semi--algebraic subsets of
$\Cal L(X)$ with values in $\hat \Cal M$ such that $[X_{n,D,W}]/L^{(n+1)d}$ is
just the volume of $Y_{n,D,W}$.  See also \S 5 and the appendix.
\bigskip
\noindent
{\bf 1.8. Remark.}  As far as we know it is not clear whether or not the
natural
morphism $\Cal M_L \rightarrow \hat \Cal M$ is injective; its kernel is
$\cap_{m \in \Bbb Z} F_m$.  However for an algebraic variety $V$ we have that
$H(V)$ and $\chi(V)$ only depend on the image of $[V]$ in $\hat \Cal M$, see
1.12.
\bigskip
\noindent
\proclaim
{1.9.  Theorem {\rm [Kon][DL3, (6.5)]}}  Using the notation of 1.4 we have the
following formula for $\Cal E_W (D)$ in terms of a log resolution $h : Y
\rightarrow X$ of supp $D$ :
$$\Cal E_W(D) = L^{-d} \sum_{I \subset T} [E^\circ_I \cap h^{-1}W] \prod_{i \in
I} \frac{L-1}{L^{\nu_i + N_i} - 1} \text { in } \hat \Cal M.$$
In particular $\Cal E_W(D)$ belongs to the image of $\Cal M_L [(L^{n}-1)^{-1}]_{n \in \Bbb N \setminus \{ 0 \}}$ in $\hat \Cal M$.
\endproclaim
\bigskip
\noindent
So by Theorem 1.4 we obtain that $\Cal E_W(D) = \Cal Z_W (D) |_{s=1}$ in $\hat
\Cal M$, where the evaluation `$s=1$' means substituting $L^{-1}$ for the
variable $L^{-s}$.
\bigskip
\noindent
{\bf 1.10.}  The following important change of variables formula is a special
case of
[DL3, Lemma 3.3], and was also mentioned in [Kon].
\bigskip
\noindent
\proclaim
{Theorem}  Let also $X^\prime$ be a smooth irreducible variety and $\rho :
X^\prime
\rightarrow X$ a proper birational morphism.  Let $D$ be an effective divisor
on $X$.  Then
$$\Cal E_W (X,D) = \Cal E_{\rho^{-1}W} (X^\prime, \rho^\ast D +
K_{X^\prime|X})$$
where $K_{X^\prime|X} = K_{X^\prime} - \rho^\ast K_X$ is the relative canonical
divisor or
{\it discrepancy divisor}.
\endproclaim
\bigskip
\noindent
{\bf 1.11.}  It is possible to generalize the set--up in 1.7 -- 1.10 to effective
$\Bbb
Q$--divisors.  We treat this in the appendix.  In particular we obtain for
an effective $\Bbb Q$--divisor $D$ on $X$, such that $rD$ is a divisor for an $r
\in \Bbb N \setminus \{ 0 \}$, an analogous invariant $\Cal E_W(D) \in \hat
\Cal M [L^{1/r}]$.  It is given in terms of a log resolution $h : Y \rightarrow
X$ (as in 1.4) by the same formula as in 1.9, where now the $N_i$ belong to
$\frac 1r (\Bbb N \setminus \{ 0 \}$).  So $\Cal E_W (D)$ belongs to the image
of $\Cal M[L^{-1/r}][(L^{n/r} - 1)^{-1}]_{n \in \Bbb N \setminus \{ 0 \}}$ in
$\hat \Cal M [L^{1/r}]$. 

When $\supp D$ has strict normal crossings we extend in the appendix the notion of
$\Cal E_W (D)$ further to the case that all coefficients of $D$ are $> -1$.
Remark that then in Theorem 1.9 (with $h = Id_X$) all $\nu_i = 1$, and our
condition on the coefficients of $D$ is thus precisely that all $\nu_i + N_i >
0$.
\bigskip
\noindent
{\bf 1.12.}  One can also specialize the invariant $\Cal E_W(D)$ to the level
of
Hodge polynomials and Euler characteristics.  We only consider expressions in
terms of log resolutions (using the notation of 1.4). 

The morphism $H : \Cal M \rightarrow \Bbb Z[u,v]$ extends canonically to a
morphism $H : \Cal M_L[(L^n-1)^{-1}]_{n \in \Bbb N \setminus \{ 0 \}}
\rightarrow \Bbb Z[uv]_{uv} [((uv)^n - 1)^{-1}]_{n \in \Bbb N \setminus \{ 0
\}} \subset \Bbb Q(u,v)$.  Since the kernel of the natural map $\Cal M_L
\rightarrow \hat \Cal M$ is killed by $H$ we can in fact consider $H$ as a
morphism from the image of $\Cal M_L[(L^n  - 1)^{-1}]_{n \in \Bbb N \setminus
\{ 0 \}}$ in $\hat \Cal M$ into $\Bbb Q(u,v)$.

We define for an effective divisor $D$ on $X$ the invariants 
$$\aligned
E_W (D) = E_W(X,D) &  := H(\Cal E_W(D)) \\
& = (uv)^{-d} \sum_{I \subset T} H(E^\circ_I \cap h^{-1} W) \prod_{i \in I}
\frac{uv-1}{(uv)^{\nu_i + N_i} - 1} \in \Bbb Q(u,v)
\endaligned
$$
and
$$e_W(D) = e_W (X,D) := \lim_{u,v \rightarrow 1} E_W(X,D) = \sum_{I \subset T}
\chi (E^\circ_I \cap h^{-1} W) \prod_{i \in I} \frac{1}{\nu_i + N_i} \in \Bbb
Q.$$ 

The extended notions of $\Cal E_W(D)$ for $\Bbb Q$--divisors of 1.11 can analogously be specialized.  We obtain the same expressions where now the $N_i$ are
rational; then $E_W(D)$ is a rational function in $u,v$ with `fractional
powers'.  For $W=X$ this was already considered by Batyrev [B3].
\bigskip
\bigskip
\head
2. Immediate generalizations and applications
\endhead
\bigskip
\noindent
{\bf 2.1.}  We recall some terminology with origins in the Minimal Model Program. See for example [KM, KMM, Kol].

On any normal variety $V$ there is a well--defined linear equivalence class of
canonical Weil divisors, denoted by $K_V$.  An arbitrary Weil divisor $D$ on
$V$ is called $\Bbb Q$--Cartier if $rD$ is Cartier for some $r \in \Bbb N
\setminus \{ 0 \}$.  A normal variety $V$ is called ($\Bbb Q$--)Gorenstein if
$K_V$ is ($\Bbb Q$--)Cartier.  

Let $X$ be a normal variety and $D$ a $\Bbb Q$--divisor on $X$ such that $K_X +
D$ is $\Bbb Q$--Cartier.  (In particular we can have $D = 0$ and then $X$ is
$\Bbb Q$--Gorenstein.)  Let $\rho : Y \rightarrow X$ be a log resolution of $\supp
D$ and denote by $E_i, i \in T$, the irreducible components of \break
$h^{-1}$ $(X_{\sing} \cup \supp D)$.  Then we can write 
$$K_Y = \rho^\ast (K_X + D) + \sum_{i \in T} (a_i - 1)E_i$$
in $\Pic Y \otimes \Bbb Q$ and $a_i = a_i (X, D;E_i)$ is called the {\it log
discrepancy} (with respect to the pair $(X,D)$) of $E_i$ for $i \in T$.  This
number $a_i$ does not depend on the chosen resolution (it is determined by the
valuation on
$k(X)$ associated to $E_i$).  Remark that when $X$ is smooth and $D=0$ the
numbers $\nu_i$ defined in 1.4 are just log discrepancies.
\medskip
($i$) Let first $D=0$.  The variety $X$ is called {\it terminal}, {\it
canonical},
{\it log
terminal} and {\it log canonical} if for some (or, equivalently, any) log
resolution
of $X$ we have
that $a_i > 1$,  $a_i \geq 1$, $a_i > 0$ and $a_i \geq 0$, respectively, for
all $i \in T$.

($ii$) When $D \ne 0$ the pair $(X,D)$ is said to be {\it Kawamata log terminal}
(shortly {\it klt}) if for some (or any) log resolution of supp $D$
we have that $a_i > 0$ for all $i \in T$.  In particular this implies that, if
$D = \sum_i d_i D_i$ with the $D_i$ irreducible, all $d_i < 1$.  (See [Kol, S]
for a discussion of other log terminality notions for pairs.)

($iii$) A closed subvariety $C \subset X$ is called a {\it log canonical centre}
of $X$ if for some log resolution $\rho : Y \rightarrow X$ there exists $i \in
T$ such that $\rho(E_i) = C$ and $a_i \leq 0$.  The {\it locus of log canonical
singularities} of $X$, denoted by $LCS(X)$, is the union of all log canonical
centres of $X$. In particular $LCS(X) = \emptyset \Leftrightarrow X$ is log
terminal.
(Hence a more appropriate notation for this locus, proposed by Koll\'ar,
would be Nlt(X), indicating the locus where $X$ is not log terminal.)
\bigskip
\noindent
{\bf 2.2.}  A natural idea, inspired by Theorem 1.10, to generalize the
invariant
$\Cal E_W(X,D)$ to a ($\Bbb Q$--)divisor $D$ on a singular variety $X$ is as
follows.  Take a resolution $h : Y \rightarrow X$ of $X$ and define
$\Cal E_W(X,D)$ as $\Cal E_{h^{-1}W} (Y,h^\ast D + K_{Y|X})$, whenever this
makes
sense, and verify independency of the chosen resolution.  So we want $X$ to be
$\Bbb Q$--Gorenstein and $h^\ast D + K_{Y|X}$ to be effective, or at least that
its coefficients are $> -1$ if its support has normal crossings. 

Below we treat the `instructional' case that $X$ is ($\Bbb Q$--)Gorenstein and
canonical
and $D$ is an effective ($\Bbb Q$--)Cartier divisor.
\bigskip
\noindent
\proclaim
{2.3. Definition -- Proposition} (i) Let $X$ be a Gorenstein and canonical
variety and $W$ a subvariety of $X$; let $D$ be an effective Cartier divisor on
$X$.  Take a resolution $h : Y \rightarrow X$ of $X$.  Then we define
$$\Cal E_W (X,D) := \Cal E_{h^{-1}W} (Y, h^\ast D + K_{Y|X}) \in {\Cal {\hat
M}}.$$

(ii) More generally let $X$ be $\Bbb Q$--Gorenstein and canonical and $W$ a
subvariety of $X$; let $D$ be an effective $\Bbb Q$--Cartier divisor on $X$.  Say $rK_X$
and $rD$ are Cartier for an $r \in \Bbb N \setminus \{ 0 \}$.  Take a
resolution $h : Y \rightarrow X$ of $X$.  Then we define $\Cal E_W (X,D) \in
\hat \Cal M [L^{1/r}]$ as above. 
\endproclaim
\bigskip
\demo{Proof} (i) The divisor $h^\ast D + K_{Y|X}$ is effective since $K_{Y|X}$
is
effective, which is equivalent to $X$ being canonical.  Let now $h^\prime :
Y^\prime \rightarrow X$ be another log resolution of $X$.
Since two such resolutions are always dominated by a third it is sufficient to
consider the case that $h^\prime$ factors through $h$ as $h^\prime : Y^\prime
\overset \pi \to \longrightarrow Y \overset h \to \longrightarrow X$.
Then by Theorem 1.10 we have
$$\aligned
\Cal E_{h^{-1}W} (Y,h^\ast D + K_{Y|X})  & = \Cal E_{\pi^{-1}h^{-1}W}
(Y^\prime,
\pi^\ast (h^\ast D + K_{Y|X}) + K_{Y^\prime|Y}) \\
& = \Cal E_{h^{\prime^{-1}}W} (Y^\prime, h^{\prime \ast} D + K_{Y^\prime|X}).
\endaligned
$$

(ii) Completely analogous, using the extended theory for $\Bbb Q$--divisors
mentioned in 1.11. \qed
\enddemo
\bigskip
When $h : Y \rightarrow X$ is a log resolution of $\supp D$ we have
the same formula as in Theorem 1.9, where the $\nu_i$ must be generalized
according to their meaning as log discrepancies.  More precisely, denoting the
irreducible components of $h^{-1}(X_{\sing} \cup \supp D)$ by $E_i, i \in T$,
we set $h^\ast D = \sum_{i \in T} N_i E_i$ and $K_Y = h^\ast K_X + \sum_{i \in
T} (\nu_i - 1)E_i$.  Then $h^\ast D + K_{Y|X} = \sum_{i \in T} (\nu_i + N_i -
1)E_i$ and so
$$\Cal E_W(X,D) = L^{-d} \sum_{I \subset T} [E^\circ_I \cap h^{-1}W] \prod_{i
\in I} \frac{L-1}{L^{\nu_i + N_i}-1},$$
where $d$ is the dimension of $X$.
\bigskip
\noindent
{\bf 2.4.} With essentially the same arguments, but needing more material from
the
appendix, we could introduce $\Cal E_W(X,D)$ for a $\Bbb Q$--Gorenstein variety
$X$ and a $\Bbb Q$--Cartier divisor $D$ on $X$ such that the pair $(X,-D)$ is
klt.  (Check that this is more general than the case in 2.3 !).  On the level
of Hodge polynomials this would be possible using [B2, Theorems 6.27 and 6.28].
We do not pursue this here; our invariants $E_W(X,D)$ in \S 3 and $\Cal E_W
(X,D)$ in \S 6 cover this case anyhow.   
\bigskip
\noindent
{\bf 2.5.}  In the rest of this section we present an application on minimal
models, taking $k = \Bbb C$.

Recall that an irreducible projective variety $V$ is called a {\it minimal
model} if $V$ is terminal and $K_V$ is {\it numerically effective} (shortly
{\it nef}\,), i.e. the intersection number $K_V \cdot C \geq 0$ for any
irreducible
curve $C$ on $V$.  The Minimal Model Program predicts the existence of a
minimal model in every birational equivalence class $\Cal C$ of nonnegative
Kodaira dimension; furthermore one should be able to transform every smooth irreducible
projective variety in $\Cal C$ by a finite number of divisorial contractions
and flips to a minimal model.

In dimension 2 it is well known that each such class has a unique minimal
model, which is moreover smooth (then divisorial contractions are just
blowing--downs and flips do not occur).  In dimension 3 the existence and
desired property of minimal models were proved by Mori; here it is crucial to
allow terminal singularities, and minimal models are {\it not} unique in a
given
birational equivalence class of nonnegative Kodaira dimension.  In dimension
$\geq 4$ the Minimal Model Program is still a major conjecture in algebraic
geometry and is becoming a working hypothesis.

It is natural and important in this context to look for invariants which are
shared by birationally equivalent minimal models.  In [Wa] Wang proved that
birationally equivalent {\it smooth} minimal models have the same Betti
numbers,
using the following result [Wa, Corollary 1.10].
\bigskip
\noindent
\proclaim
{2.6. Proposition} Let $f : V -\!\rightarrow \!V^\prime$ be a birational map
between two
minimal models.  Then there exist a smooth projective variety $Y$ and
birational morphisms $\varphi : Y \rightarrow V, \varphi^\prime : Y \rightarrow
V^\prime$ such that $\varphi^\ast K_V = \varphi^{\prime \ast} K_{V^\prime}$.
\endproclaim
\bigskip
\noindent
(In fact Wang only needs $V$ and $V^\prime$ to be terminal varieties for which
$K_V$ and $K_{V^\prime}$ are nef along the exceptional loci of $f$ in $V$ and
$V^\prime$, respectively, to conclude.)  This result has more interesting
consequences.
\bigskip
\noindent
\proclaim
{2.7. Theorem} Let $V$ and $V^\prime$ be birationally equivalent minimal
models.  Then

(i) $\Cal E_V(V,0) = \Cal E_{V'}(V^\prime,0)$, and

(ii) if $V$ and $V^\prime$ are smooth, then $[V] = [V^\prime]$.
\endproclaim
\bigskip
\demo{Proof} (i) Take $V \overset \varphi \to \longleftarrow Y \overset
\varphi^\prime \to \longrightarrow V^\prime$ as in Proposition 2.6.  Then by
Theorem 1.10 (and its generalization in 1.11) we have
$$\Cal E_V(V,0) = \Cal E_Y(Y,K_Y - \varphi^\ast K_V) = \Cal E_Y(Y,K_Y -
\varphi{^\prime \ast} K_{V^\prime}) = \Cal E_{V'}(V^\prime,0).$$

(ii) For any smooth variety $X$ we have that $\Cal E_X(X,0) = [X]$.
\qed
\enddemo
\bigskip
\noindent
As a corollary birationally equivalent {\it smooth} minimal models have the
same Hodge numbers and a  fortiori the same Betti numbers.  In particular this
is true for smooth Calabi--Yau varieties. See also [B1, Theorems 1.1 and 4.2].
\bigskip
\noindent
{\bf 2.8.}  Assuming the Minimal Model Program in some dimension $d$ we can use
Theorem 2.7 to define a birational invariant.  For any birational equivalence
class ${\Cal C}$ of nonnegative Kodaira dimension the expression $\Cal E
:= \Cal E_X(X,0)$ is independent of a chosen minimal model $X$.  Looking
at 2.3 it is given by the following formula in terms of any log resolution $h :
Y \rightarrow X$ of any minimal model $X$ of $\Cal C$.  Denote by $E_i, i \in
T$, the irreducible components of $h^{-1}(X_{\text{sing}})$ and set $K_Y =
h^\ast K_X + \sum_{i \in T} (\nu_i - 1)E_i$.  Then
$${\Cal E} = L^{-d} \sum_{I \subset T} [E^\circ_I]  \prod_{i \in I}
\frac{L-1}{L^{\nu_i} - 1}.$$

\smallskip
One could extract \lq minimal stringy Hodge numbers\rq\ from (the Hodge polynomial version of) this invariant, see [B2]; and maybe it is related to a \lq minimal cohomology theory\rq\ as explained in [Wa].
\bigskip
\bigskip
\head
3. Singular varieties; on the level of Hodge polynomials and Euler
characteristics
\endhead
\bigskip
\noindent
{\bf 3.1.} Our aim in this paper is to associate zeta functions and
`Kontsevich'
invariants to effective $\Bbb Q$--Cartier divisors $D$ on arbitrary $\Bbb
Q$--Gorenstein varieties $X$ for which $X_{\sing} \subset \supp D$, generalizing
the notions in \S 1.  In this section we realize this on the level of Hodge
polynomials and Euler characteristics in a fairly elementary way.  The more
general case on the level of the Grothendieck ring will be treated in \S 5.
\bigskip
\noindent
{\bf 3.2.} We fix notation for this section.  Let $X$ be a $\Bbb Q$--Gorenstein
variety and $D$ an effective $\Bbb Q$--Cartier divisor on $X$. 
(When $\dim X =2$ we only need that $X$ is normal and $D$ can be any effective Weil divisor with rational coefficients, see [V3].)
For a log
resolution $h : Y \rightarrow X$ of $\supp D$ we denote by $E_i, i
\in T$, the irreducible components of $h^{-1}(X_{\sing} \cup \supp D)$ and we
put $E^\circ_I := (\cap_{i \in I} E_i) \setminus (\cup_{\ell \not\in I}
E_\ell)$ for $I \subset T$.  We also set 
$h^\ast D = \sum_{i \in T} N_i E_i$ and $K_Y = h^\ast K_X + \sum_{i \in T}
(\nu_i - 1)E_i$.  Remember that now the $\nu_i \in \Bbb Q$ and they can be
negative or zero.

In the sequel we will again consider arbitrary subvarieties $W$ of $X$. One can think mainly about $W$ being for example $X$, $\supp D$, $X_{\sing}$ or a point of $X_{\sing}$.
\bigskip
\noindent
\proclaim
{\bf 3.3. Definition -- Proposition}  Let $X$ be a $\Bbb Q$--Gorenstein variety
of dimension $d$ and $W$ a subvariety of $X$.  Let $D$ be an effective $\Bbb
Q$--Cartier divisor on $X$ such that $X_{\sing} \subset \supp D$.  Take $r \in
\Bbb N \setminus \{ 0 \}$ with $rK_X$ and $rD$ Cartier.  

(i) The zeta function $Z_W (D,s) = Z_W (X,D,s)$ is the unique rational function
in the variable $(uv)^{-s/r}$ and with coefficients in (the fraction field of)
$\Bbb Z[u,v][(uv)^{1/r}]$ such that for $n >> 0$
$$Z_W(D,n) = E_{h^{-1}W} (Y,nh^\ast D + K_{Y|X}),$$
where $h : Y \rightarrow X$ is a resolution of $X$.

(ii) Let $h : Y \rightarrow X$ be a log resolution of $\supp D$.
With the notation of 3.2 we have that
$$Z_W (D,s) = \frac{1}{(uv)^d} \sum_{I \subset T} H(E^\circ_I \cap h^{-1}W)
\prod_{i \in I} \frac{uv - 1}{(uv)^{\nu_i+sN_i}-1} .$$ 
\endproclaim
\bigskip
\demo{Proof}  Let $h_i : Y_i \rightarrow X$ be resolutions of $X$ for $i =
1,2$.  We first show that the defining expressions for $Z_W (D,n)$ using $Y_1$
and $Y_2$ are equal when $n >> 0$.  Take $n$ such that $n h^\ast_i D +
K_{Y_i|X}$ is effective for $i = 1,2$ (here we need that $X_{\sing} \subset
\supp D)$, and take a resolution $h : Y \rightarrow X$ of $X$ dominating both
$Y_1$ and $Y_2$.
\vskip 5 truemm
\centerline{
\beginpicture
\setcoordinatesystem units <.5truecm,.5truecm>
\put{$\searrow$} at 1 1
\put{$\swarrow$} at 1 -1
\put{$\searrow$} at -1 -1
\put{$\swarrow$} at -1 1
\put{$\Bigg\downarrow$} at 0 0
\put{$\varphi_1$} at -1.4 1.4
\put{$\varphi_2$} at 1.4 1.4
\put{$h$} at -.4 0
\put{$h_2$} at 1.5 -1.4
\put{$h_1$} at -1.4 -1.4
\put{$Y_1$} at -2 0
\put{$Y$} at  0 2
\put{$Y_2$} at 2 0
\put{$X$} at  0 -2
\endpicture}
\vskip 5 truemm
\noindent
Then by Theorem 1.10 (for $\Bbb Q$--divisors and on the level of Hodge
polynomials) we have for $i = 1,2$ that
$$\aligned
E_{h^{-1}_i W} (Y_i, nh^\ast_i D + K_{Y_i|X}) & = E_{\varphi^{-1}_i h^{-1}_iW}
(Y, \varphi^\ast_i (nh^\ast_iD + K_{Y_i|X}) + K_{Y|Y_i}) \\
& = E_{h^{-1}W} (Y,nh^\ast D + K_{Y|X}).
\endaligned
$$
Choosing now $h : Y \rightarrow X$ as a log resolution for $\supp D$
we have that 
$$
E_{h^{-1}W} (Y,nh^\ast D + K_{Y|X}) = \frac{1}{(uv)^d} \sum_{I \subset T}
H(E^\circ_I \cap h^{-1}W) \prod_{i \in I} \frac{uv-1}{(uv)^{\nu_i+nN_i}-1}.$$ 
Hence for $n >> 0$ the stated rational function in (ii) indeed yields
$E_{h^{-1}W} (Y,nh^\ast D + K_{Y|X})$ when evaluating in $s=n$ (i.e. in
$(uv)^{-s/r} = (uv)^{-n/r})$.

Finally this rational function must be unique since a polynomial over the
domain $\Bbb Z[u,v][(uv)^{1/r}]$ can have at most finitely many zeroes. \qed
\enddemo
\bigskip
\proclaim
{3.4. Definition}  With the same notation as in 3.3 we define the {\it
topological zeta function} of $D$ as
$$z_W(D,s) = z_W(X,D,s) := \sum_{I \subset T} \chi (E^\circ_I \cap h^{-1} W)
\prod_{i \in I} \frac{1}{\nu_i + sN_i} \in \Bbb Q(s).$$
We can justify this definition either by an analogous proof or by obtaining
$z_W(D,s)$ from $Z_W(D,s)$ by a limit argument as in 1.6.
\endproclaim
\bigskip
\noindent
{\bf 3.5.} In the following we extend Kontsevich's construction $E_W(X,D)$ to
$\Bbb
Q$--Gorenstein varieties $X$.  We should remark here that in [DL3] Denef and
Loeser also generalized in a different way $\Cal E_W(X,D)$ to (arbitrary)
singular varieties $X$.  We consider their point of view as more
`integrational' and ours as more `geometrical'.  Our idea is simply to
substitute $s=1$ in $Z_W(D,s)$ when this makes sense or, more generally, to
take the limit for $s \rightarrow 1$.
\bigskip
\noindent
\proclaim
{3.6. Definition}  Let $X$ be a $\Bbb Q$--Gorenstein variety and $W$ a
subvariety
of $X$.  Let $D$ be an effective $\Bbb Q$--Cartier divisor on $X$  such that
$X_{\sing} \subset
\supp D$.  Take $r \in \Bbb N \setminus \{ 0 \}$ with $rK_X$ and $rD$ Cartier.
Then we put
$$E_W(D) = E_W (X,D) := \lim_{s \rightarrow 1} Z_W(X,D,s) \in \Bbb
Q(u^{1/r},v^{1/r}) \cup \{ \infty \}.$$ 
\endproclaim
\bigskip
\noindent
{\sl Remarks.} (1) By $\lim_{s \rightarrow 1}$ we mean taking the limit
$(uv)^{-s/r}
\rightarrow (uv)^{-1/r}$. This is well defined since $Z_W(D,s)$ is a rational function in the
variable $(uv)^{-s/r}$ over a field.

(2) If there exists a log resolution $h : Y \rightarrow X$ of $\supp D$ for
which $\nu_i + N_i \ne 0$ for all $i \in T$, then, because of the formula in
3.3(ii), we obtain $E_W(D)$ from $Z_W(D,s)$ simply by {\it substituting}
$(uv)^{-1/r}$ for $(uv)^{-s/r}$.  (We formulate this below as Proposition 3.7.)
If on the other hand there does not exist such a log resolution, then in
general we will have $E_W(D) = \infty$.  However there are cases where our
definition then yields an element in $\Bbb Q(u^{1/r},v^{1/r})$, see example 4.1.
\bigskip
\noindent
\proclaim
{3.7. Proposition}  Let $W \subset X$ and $D$ be as in 3.6.  Let $h : Y
\rightarrow X$ be a log resolution of $\supp D$ for which $\nu_i +
N_i \ne 0$ for all $i \in T$ (using the notation of 3.2).  Then
$$E_W(D) = \frac{1}{(uv)^d} \sum_{I \subset T} H(E^\circ_I \cap h^{-1}W)
\prod_{i \in I} \frac{uv-1}{(uv)^{\nu_i+N_i}-1}.$$
\endproclaim

\noindent
So indeed we extended Kontsevich's invariant for smooth $X$ on the level of
Hodge polynomials (1.12).
\bigskip
\noindent
\proclaim
{3.8. Definition -- Proposition}  Let $W \subset X$ and $D$ be as in 3.6.  We
define
$$e_W(D) = e_W(X,D) := \lim_{s \rightarrow 1} z_W(X,D,s) \in \Bbb Q \cup \{
\infty \}.$$
Let $h : Y \rightarrow X$ be a log resolution of $\supp D$ for which
$\nu_i + N_i \ne 0$ for all $i \in T$.  Then
$$e_W(D) = \sum_{I \subset T} \chi (E^\circ_I \cap h^{-1}W) \prod_{i \in I}
\frac{1}{\nu_i + N_i}.$$
\endproclaim
\bigskip
\noindent
{\bf 3.9.} Next we introduce analogous invariants for {\it pairs} $(X,D)$,
which
will coincide with Batyrev's {\it stringy $E$--function} and {\it stringy  Euler number} for klt pairs
[B3].
\bigskip
\noindent
\proclaim
{3.10. Definition -- Proposition}  Let $X$ be a $\Bbb Q$--Gorenstein variety and
$W$ a subvariety of $X$.  Let $D$ be an effective $\Bbb Q$--Cartier divisor on
$X$ such that
$X_{\sing}
\subset \supp D$.  Take $r \in \Bbb N \setminus \{ 0 \}$ with $rK_X$ and $rD$
Cartier.

(i) We put
$$E_W \bigl((X,D)\bigr) := \lim_{s \rightarrow -1} Z_W(X,D,s) \in \Bbb Q
(u^{1/r},v^{1/r}) \cup \{ \infty \}.$$

(ii) Let $h : Y \rightarrow X$ be a log resolution of $\supp D$.
Using the notation of 3.2, let $a_i, i \in T$, denote the log discrepancy of
$E_i$ with respect to the pair $(X,D)$. Then, if $a_i \ne 0$ for all $i \in T$,
we have
$$E_W \bigl((X,D)\bigr) = \frac{1}{(uv)^d} \sum_{I \subset T} H(E^\circ_I \cap h^{-1}W)
\prod_{i \in I}
\frac{uv-1}{(uv)^{a_i}-1}.$$
\endproclaim
\bigskip
\noindent
{\sl Remark.}  By $\lim_{s \rightarrow -1}$ we mean taking the limit
$(uv)^{-s/r}
\rightarrow (uv)^{1/r}$.
\bigskip
\noindent
\demo {Proof}  If $\nu_i - N_i \ne 0$ for all $i \in T$, then, because of the
formula
for $Z_W(X,D,s)$ in 3.3(ii) this limit procedure just means substituting
$(uv)^{1/r}$ for the variable $(uv)^{-s/r}$.  Clearly we obtain the stated
formula for $E_W \bigl((X,D)\bigr)$ since $a_i  = \nu_i - N_i$ for $i \in T$. \qed
\enddemo
\bigskip
\noindent
{\bf 3.11.}  When the pair $(X,D)$ is 
klt and for $W=X$ Batyrev introduced in [B3] the
same invariant as the {\it stringy $E$--function} of $(X,D)$, denoted by
$E_{st}(X,D)$.  (We do not recover his invariant completely as a special case
of $E_W \bigl((X,D)\bigr)$ because Batyrev only requires $K_X + D$ to be $\Bbb
Q$--Cartier.)  Analogously the invariant $e_X \bigl((X,D)\bigr)$ below was baptized {\it
stringy Euler number} by Batyrev and denoted by $e_{st}(X,D)$.
\bigskip
\noindent
\proclaim
{3.12. Definition -- Proposition}  Let $W \subset X$ and $D$ be as in 3.10.  We
define
$$e_W \bigl((X,D)\bigr) := \lim_{s \rightarrow -1} z_W (X,D,s) \in \Bbb Q \cup \{ \infty
\}.$$
Let $h : Y \rightarrow X$ be a log resolution of $\supp D$ for which
$\nu_i - N_i \ne 0$ for all $i \in T$.  Then, denoting by $a_i$ the log
discrepancy of $E_i$ with respect to the pair $(X,D)$, we have
$$e_W \bigl((X,D)\bigr) = \sum_{I \subset T} \chi (E^\circ_I \cap h^{-1}W) \prod_{i \in I}
\frac{1}{a_i}.$$
\endproclaim
\bigskip
\noindent
{\bf 3.13.}  In Definition--Proposition 3.3, and hence in all subsequent
constructions, we required the effective divisor $D$ to satisfy $X_{\sing}
\subset \supp D$.  We needed this to assure that for a resolution $h : Y
\rightarrow X$ the divisor $nh^\ast D + K_{Y|X}$ would be effective for $n >>
0$.  However, using $A5$ in the appendix, it is in fact sufficient to require that $\supp D$
contains the locus of log canonical singularities $LCS(X)$ of $X$. 
\bigskip
\noindent
\proclaim
{3.14. Definition -- Theorem}  Let $X$ be a $\Bbb Q$--Gorenstein variety of
dimension $d$ and $W$ a subvariety of $X$.  Let $D$ be an effective $\Bbb
Q$--Cartier divisor on $X$ such that 
$LCS (X) \subset \supp D$.  Take $r \in \Bbb N \setminus \{  0 \}$ with $rK_X$
and $rD$ Cartier.

(i) The zeta function $Z_W(D,s) = Z_W (X,D,s)$ is the unique rational function
in the variable $(uv)^{-s/r}$ and with coefficients in (the fraction field of)
$\Bbb Z[u,v][(uv)^{1/r}]$ such that for $n >> 0$
$$Z_W(D,n) = E_{h^{-1}W} (Y,nh^\ast D + K_{Y|X}),$$
where $h : Y \rightarrow X$ is a {\it log resolution} for $\supp D$.

(ii) With the notation of 3.2 for $h$ we have that
$$Z_W(D,s) = \frac{1}{(uv)^d} \sum_{I \subset T} H(E^\circ_I \cap h^{-1}W)
\prod_{i \in T} \frac{uv-1}{(uv)^{\nu_i + sN_i} - 1}.$$
\endproclaim
\bigskip
\noindent
\demo{Proof}  We proceed analogously as in the proof of 3.3, but now working
only
with log resolutions $h : Y \rightarrow X$ of $\supp D$.
Then for $n >> 0$ the coefficients $d_i = nN_i + \nu_i - 1$ of $nh^\ast D +
K_{Y|X}$ all satisfy $d_i > -1$. Indeed any exceptional component $E_i$ of $h$
for which $\nu_i \leq 0$ satisfies $h(E_i) \subset LCS(X) \subset \supp D$, and
hence $N_i > 0$ for such an $E_i$.  So in this case the invariant $E_{h^{-1}W}
(Y,nh^\ast D + K_{Y|X})$ is well defined by $A5$ and Theorem A6. \qed
\enddemo
\bigskip
\noindent
{\sl Remark.} In the formula above the `denominators' $\nu_i + sN_i$ are thus
always
nonzero since either $\nu_i > 0$ or $N_i > 0$.
\bigskip
\noindent
{\bf 3.15.}  We can also extend all invariants which we considered in 3.4 -- 3.12,
i.e.
$z_W(D,s)$, $E_W(D)$, $e_W(D)$, $E_W \bigl((X,D)\bigr)$ and $e_W \bigl((X,D)\bigr)$, to the case that
only
$LCS(X) \subset \supp D$.

In particular when $X$ is log terminal and $D=0$, then our
invariants
$E_X(0)$ and $e_X(0)$ are precisely the stringy E--function $E_{st}(X;u,v)$ and
stringy Euler number $e_{st}(X)$ of Batyrev [B2]. 
\bigskip
\bigskip

\head 4. Examples
\endhead 

\bigskip
In this section we present a number of examples, first in dimension two and then in higher dimension, for which we compute the invariants introduced above. Recall (see (0.4(a)) that in dimension two we can consider more generally Weil divisors instead of Cartier divisors.
\bigskip
\noindent
{\bf 4.1.} Let $0 \in X$ be a normal surface germ with minimal resolution $h :
Y
\rightarrow X$ such that $h^{-1} \{ 0 \} = E_0 \cup E_g$, where $E_0$ and $E_g$
are nonsingular curves of genus $0$ and $g \geq 2$,  \linebreak  
respectively, intersecting
transversely.  So $h$ is
already a log resolution of $X$.  (This singularity is quasihomogeneous.)  Let
$E$ be
a nonsingular curve (germ) in $Y$ intersecting $E_g$ transversely in one point
and disjoint from $E_0$.  Denote $D = h(E)$; so $D$ is a prime Weil divisor on
$X$ through $0$ and $h$ is also a log resolution of $D$. See Figure 1.

\vskip 6truemm
\centerline{
\beginpicture
\setcoordinatesystem units <.47truecm,.47truecm>
\putrectangle corners at 0 0 and 7 6
\ellipticalarc  axes ratio 4:2  90 degrees from 2 3 center at 2 5
\ellipticalarc  axes ratio 4:2  90 degrees from 6 1 center at 2 1
\put {$\bullet$} at 2 3
\put {$D$} at 5.1 1
\put {$0$} at 1.9 3.7
\put {$X$} at 8 5
\put {$\longrightarrow$} at -3 3
\put {$h$} at -3 3.7
\setcoordinatesystem units <.47truecm,.47truecm> point at 14 0
\putrectangle corners at 0 0 and 8 6
\putrule from 2 .8 to 2 5.2
\putrule from 6 .8 to 6 5.2
\putrule from 1 4.5 to 7 4.5
\put {$E_0$} at 1.2 2
\put {$E_g$} at 4 3.7
\put {$E$} at 6.8 2
\put {$Y$} at -1 5
\put {Figure 1} at 11 -1.5
\endpicture}
\vskip 6truemm

\noindent
Let $-\kappa_0$ and $-\kappa_g$ denote the self--intersection number of
$E_0$ and $E_g$ on $Y$, respectively; we have that $\kappa_0 \geq 2$ and
$\kappa_g \geq 1$.  We will treat the germs $0 \in X$ for which $N := 2g -
\kappa_g - 1 > 0$ in order to compute $z_0 (ND,s)$ and $e_0(ND)$ for the
effective Weil divisor $ND$ on $X$.

We denote as usual $h^\ast ND = NE + N_0 E_0 + N_g E_g$ and $K_Y =
h^\ast K_X + (\nu_0 - 1)E_0 +$ $(\nu_g - 1)E_g$.  The following relations
are well known (see for example [V3, Lemma 2.3]) :
$$\left\{ \aligned
\kappa_0 N_0 & = N_g \\
\kappa_0 \nu_0 & = \nu_g + 1
\endaligned \right.
\qquad \text{ and } \qquad
\left\{ \aligned
\kappa_g N_g & = N_0 + N \\
\kappa_g \nu_g & = (\nu_0 - 1) + 2 - 2g.
\endaligned \right.
$$
A short computation yields the expression for $N_0$, $\nu_0$, $N_g$ and $\nu_g$
in
terms of our data $\kappa_0$, $\kappa_g$ and $g$ :
$$\left\{ \aligned  N_0 & = \frac{N}{\kappa_0 \kappa_g - 1} =
\frac{2g-\kappa_g-1}{\kappa_0 \kappa_g - 1} \\
  \nu_0 & = \frac{-2g + \kappa_g + 1}{\kappa_0 \kappa_g - 1}
 \endaligned \right.
 \qquad \text { and } \qquad
 \left\{ \aligned
N_g & = \frac{\kappa_0 N}{\kappa_0 \kappa_g - 1} =
\frac{\kappa_0(2g-\kappa_g - 1)}{\kappa_0 \kappa_g - 1} \\
\nu_g & = \frac{\kappa_0(1 - 2g) + 1}{\kappa_0 \kappa_g - 1}
\endaligned \right.
$$
Remark that $\nu_0 + N_0 = 0$ (which, as you can guess, is forced by our choice of $N$); nevertheless $e_0(ND)$ will be a rational  number.  We have by definition that $$\aligned z_0 (ND,s) & = \frac{1}{\nu_0 + sN_0} + \frac{1}{(\nu_0 +
sN_0)(\nu_g + sN_g)} + \frac{-2g}{\nu_g + sN_g} + \frac{1}{(\nu_g + sN_g)(1 +
sN)} \\
& = \frac{1 + (\kappa_0 - 2g)(1 + sN)}{(\nu_g + sN_g)(1 + sN)}.
\endaligned
$$
The fact that $\nu_0 + sN_0$ cancels in the denominator is a general fact; see
[V3, 2.2].  Plugging in the expression for $\nu_g$ and $N_g$ yields
$$z_0 (ND,s) = \frac{(\kappa_0 \kappa_g - 1)[1 + (\kappa_0 -
2g)(1 + sN)]}{(1 - 2 \kappa_0 g + \kappa_0 (1 + sN))(1 + sN)}
\qquad\qquad \text{
(with $N = 2g - \kappa_g - 1)$}
$$
and
$$e_0 (ND) = \lim_{s \rightarrow 1} z_0 (ND,s) = \frac{(2g -
\kappa_0)(2g - \kappa_g) - 1}{2g - \kappa_g}.$$
One can analogously compute $Z_0 (ND,s)$ and $E_0(ND)$.
\bigskip
\noindent
{\bf 4.2.} Let $0 \in X$ and $h : Y \rightarrow X$ be as above with $g = 1$
(instead
of $g \geq 2$).  Now let $E^\prime$ be a nonsingular curve germ in $Y$
intersecting $E_0$ transversely in one point and disjoint from $E_1$, and
denote $D^\prime = h(E^\prime)$. See Figure 2.

\vskip 6truemm
\centerline{
\beginpicture
\setcoordinatesystem units <.47truecm,.47truecm>
\putrectangle corners at 0 0 and 7 6
\ellipticalarc  axes ratio 4:2  90 degrees from 2 3 center at 2 5
\ellipticalarc  axes ratio 4:2  90 degrees from 6 1 center at 2 1
\put {$\bullet$} at 2 3
\put {$D'$} at 5 1
\put {$0$} at 1.9 3.7
\put {$X$} at 8 5
\put {$\longrightarrow$} at -3 3
\put {$h$} at -3 3.7
\put {Figure 2} at -3 -1.5
\setcoordinatesystem units <.47truecm,.47truecm> point at 13 0
\putrectangle corners at 0 0 and 7 6
\putrule from 2 .8 to 2 5.2
\putrule from 1 1.5 to 5 1.5
\putrule from 1 4.5 to 5 4.5
\put {$E_0$} at 1.2 3
\put {$E_1$} at 5.8 4.5
\put {$E'$} at 5.8 1.5
\put {$Y$} at -1 5
\endpicture}
\vskip 6truemm

\noindent
One easily computes (see [V3, 2.5]) that
$$z_0 (ND^\prime, s) = - \frac{\kappa_0 \kappa_1 - 1}{1 + sN}
\qquad
\text{and thus} \qquad e_0 (ND^\prime) = - \frac{\kappa_0 \kappa_1 -
1}{1
+ N}.$$

Now choose $N = \kappa_0 - 1$. It is easy to verify that then $\nu_1 + N_1
= 0$; so as in 4.1 we could not have defined $e_0 ((\kappa_0 - 1)D')$ by the usual formula. However in this example our definition on the level of Hodge
polynomials yields $E_0 ((\kappa_0 - 1)D') = \infty$.

\bigskip
\noindent
{\bf 4.2.1.} {\sl Remark.} One could argue whether in Definition--Proposition 3.8 (and analogously in 3.12) it is more appropriate to introduce $e_W(D)$ as $\lim_{u,v \rightarrow 1} E_W(D)$. When $E_W(D) \neq \infty$ this amounts to the same, but when $E_W(D) = \infty$ we then would miss some interesting values of $e_W(D)$ as in 4.2 above. 

\bigskip
\noindent
{\bf 4.3.} Let $X$ be the quadric hypersurface $\{ xy - zw = 0 \}$ in $\Bbb
A^4$.  The origin $0$ is the only singular point of $X$.  Blowing up $0$ yields
a log resolution $h_1 : Y_1 \rightarrow X$ of $X$, which is an isomorphism
outside $h^{-1} \{ 0 \}$ and with $E_1 = h^{-1} \{ 0 \} \cong ( \{ xy - zw = 0
\} \subset \Bbb P^3) \cong \Bbb P^1 \times \Bbb P^1$.
\bigskip
(a) Consider the divisor $D = E + E^\prime$ on $X$, where $E$ and $E^\prime$
are the zero sets of the functions $z - w$ and $y$ on $X$, respectively.
Remark that $E$ is irreducible and that $E^\prime$ consists of two irreducible
components.  We want to compute $z_0(D,s)$.  In this example we will use the
same notation for divisors and their strict transforms by blowing--ups.

\vskip 6truemm
\centerline{
\beginpicture
\setcoordinatesystem units <.47truecm,.47truecm>
\putrule from 0 0 to 0 8
\putrule from 4 0 to 4 8
\putrule from 5.5 2 to 9.5 2
\putrule from 9.5 2 to 9.5 10
\putrule from 11 4 to 11 12
\putrule from 1.5 8.55 to 1.5 10
\putrule from 7 10.55 to 7 12
\putrule from 0 0 to 8 0
\putrule from 1.5 10 to 9.5 10
\setlinear \plot  0 0  4 1.45 /
     \plot 4 0  5.5 2 /
	 \plot 8 0  11 4  9.5 3.45 /
	 \plot 0 8  11 12 /
	 \plot 4 8  7 12 /
\setdots <2pt>
\putrule from 1.5 2 to 1.5 8.55
\putrule from 7 4 to 7 10.55
\putrule from 1.5 2 to 5.5 2
\putrule from 3 4 to 11 4
\setlinear \plot 4 1.45  9.5 3.45 /
       \plot  0 0  3 4 /
	   \plot 5.5 2  7 4 /
\setsolid
\linethickness=.8pt
\putrule from 5.5 2 to 5.5 10
\multiput {$E'$} at 2 10.7  6 12 /
\put {$E$} at 11.7 11.5
\put {$E_1$} at 9 0
\put {Figure 3} at 6 -2
\endpicture}
\vskip 6truemm

In Figure 3 we present the intersection configuration of $E_1,E$ and $E^\prime$
on $Y_1$.  The variety $Y_1$ is naturally covered by 4 affine charts, each
isomorphic to $\Bbb A^3$.  In the `main chart' the exceptional surface $E_1$
and the strict transforms $E$ and $E^\prime$ are given in affine coordinates
$x,z,w$ by
$$\align
& E_1 :  x = 0 , \\
& E :  z - w = 0, \\
& E^\prime :  z \cdot w = 0
\endalign
$$
(in the other charts $E$ and $E^\prime$ do not intersect).

We obtain a log resolution $h$ of $D$ by composing $h_1$ with the blowing--up
$h_2 : Y_2 \rightarrow Y_1$ of the curve $E \cap E^\prime (\cong \Bbb A^1)$ in
$Y_1$.  The exceptional variety $E_2$ of $h_2$ is isomorphic to $\Bbb A^1
\times \Bbb P^1$; the intersection configuration of $E_2, E_1, E$ and
$E^\prime$ is presented in Figure 4.

\vskip 5truemm
\centerline{
\beginpicture
\setcoordinatesystem units <.47truecm,.47truecm>
\putrule from 2 2 to 2 8
\putrule from 5 4 to 5 10
\putrule from 7 2 to 7 8
\putrule from 8 6 to 8 12
\putrule from 10 4 to 10 10
\putrule from 11 8 to 11 14 
\putrule from 13 6 to 13 12
\putrule from 14 10 to 14 16
\putrule from 3 8.66 to 3 12
\putrule from 6 10.66 to 6 14
\putrule from 9 12.66 to 9 16
\setlinear 
     \plot 4 0  16 8 /
	 \plot 2 8  14 16 /
	 \plot 0 4  2 5.33 /
	 \plot 2 2  5 4  7 2 /
	 \plot 7 5.33  8 6  10 4   /
	 \plot 10 7.33  11 8  13 6  /
	 \plot 13 9.33  14 10  16 8 /
	 \plot 7 8  3 12 /
	 \plot 10 10  6 14 /
	 \plot 13 12  9 16 /
	 \plot 4 0  0 4 /
\setdots <2pt>
\putrule from 3 6 to 3 8.66
\putrule from 6 8 to 6 10.66
\putrule from 9 10 to 9 12.66
\setlinear \plot 2 5.33  12 12  14 10 /
       \plot  3 6  5 4  7 5.33 /
	   \plot 6 8  8 6  10 7.33 /
	   \plot 9 10  11 8  13 9.33 /
\setsolid
\multiput {$E'$} at 2.2 11.5  5.2 13.5 /
\put {$E$} at 8.2 15.5
\put {$E_1$} at  15.5 6.5
\put {$E_2$} at 14.7 15
\put {Figure 4} at 8 -2
\endpicture}
\vskip 6truemm

\noindent
Denoting as usual $h^\ast D =  E + E^\prime + N_1 E_1 + N_2 E_2$ and $K_Y =
h^\ast K_X + (\nu_1 - 1)E_1 + (\nu_2 - 1)E_2$, one easily verifies that
$(\nu_1,N_1) = (2,2)$ and $(\nu_2,N_2) = (2,3)$.  The contributors to
$z_0(D,s)$ are $E^\circ_1,(E_1 \cap E_2)^\circ, (E_1 \cap E)^\circ, (E_1 \cap
E^\prime)^\circ, E_1
\cap E_2
\cap E$ and $E_1 \cap E_2 \cap E^\prime$.  Now $\chi(E^\circ_1) = 0$ and the
other Euler characteristics are obvious; then
$$\align
z_0 (D,s) & = \frac{1}{\nu_1 + sN_1} \big(\frac{-1}{\nu_2 + sN_2} + \frac{3}{1+s} +
\frac {3}{(\nu_2 + sN_2)(1+s)}\big) \\
& = \frac{4}{(2+3s)(1+s)} .
\endalign
$$
Also $e_0(D) = \lim_{s \rightarrow 1} z_0(D,s) = \frac 25$ and $e_0\big((X,D)\big) =
\lim_{s \rightarrow -1} z_0(D,s) = \infty.$
\bigskip
(b) Now consider the $\Bbb Q$--divisor $D = NE + N^\prime E^\prime$ with $N > 0,
N^\prime > 0, N \ne 1, N^\prime \ne 1$ and $N + N^\prime = 2$.  The morphism $h
:
Y_2 \rightarrow X$ in (a) is of course  still a log resolution of $D$.  The
only difference with the data in (a) is that here $h^\ast D = NE + NE^\prime +
N_1 E_1 + N_2 E_2$ with $N_1 = N + N^\prime = 2$ and $N_2 = N + 2N^\prime = 2 +
N^\prime$.  So
$$\align
z_0(D,s)  = & \frac{1}{\nu_1 + sN_1} \big(\frac{-1}{\nu_2 + sN_2} + \frac{1}{1 +
sN}
+ \frac{2}{1 + sN^\prime} + \frac{1}{(\nu_2 + sN_2)(1
+ sN)} \\ & \qquad \qquad \qquad \qquad \qquad \qquad \quad + \frac {2}{(\nu_2
+
sN_2)(1 +
sN^\prime)}\big) \\
 = & \frac{1}{2 + 2s} \cdot \frac{8 + 16s + 8s^2}{(2 + s(2 + N^\prime))(1 +
 sN)(1
+ sN^\prime)} \\
 = & \frac{4(1+s)}{(2 + s(2 + N^\prime))(1+sN)(1+sN^\prime)} \, .
\endalign
$$
And then $e_0(D,s) = \frac{8}{(4 + N^\prime)(1+N)(1+N^\prime)}$ and $e_0\big((X,D)\big)
= 0$.

\bigskip
\noindent
{\bf 4.4.} Fix $d \in \Bbb N, d \geq 3$.  Take a homogeneous polynomial $F$ in
$d+1$
variables of degree $a  \geq 2$ such that $\{ F = 0 \} \subset \Bbb P^d$ is
nonsingular.

Let $X$ be the hypersurface in $\Bbb A^{d+1}$ given by the zero set of $F$; so
$X$ is the affine cone over $\{ F = 0 \} \subset \Bbb P^d$ and the origin is
the only singular point of $X$.  Let $D$ be the intersection of $X$ with a
general hyperplane through the origin in $\Bbb A^{d+1}$.  The blowing--up $h : Y
\rightarrow X$ of the origin yields a log resolution of $X$, which is moreover
a log resolution of $D$.  We denote the strict transform of $D$ by $E$, and the
exceptional variety of $h$ by $E_1$.  Notice that $E_1$ is isomorphic to $\{ F
= 0 \} \subset \Bbb P^d$.  We try to give an impression of this situation in
Figure 5. 

\vskip 6truemm
\centerline{
\beginpicture
\setcoordinatesystem units <.47truecm,.47truecm>
\setlinear \plot  -1 3  1 5  -1 -5  1 -3  -1 3 /
\ellipticalarc axes ratio 3:1.1  360 degrees from 3 4 center at 0 4
\ellipticalarc axes ratio 3:1.1  360 degrees from 3 -4 center at 0 -4
\setquadratic \plot  3 4  2.2 2  0 0  -2.2 -2  -3 -4 /
        \plot  -3 4  -2.2 2  0 0  2.2 -2  3 -4 /
\put {$\bullet$} at 0 0
\put {$D$} at 1 2
\put {$0$} at -1 0
\put {$X$} at -3.6 3
\setlinear
\setshadegrid span <2truept>
\hshade -5 -1 -1   -3 -.6 1   0 0 0   3 -1 .6   5 1 1  /   
\put {$\longleftarrow$} at 8 0
\put {$h$} at 8 .7
\put {Figure 5} at 8 -6.5
\setcoordinatesystem units <.47truecm,.47truecm> point at -16 0
\ellipticalarc axes ratio 3:1.1  360 degrees from 3 4 center at 0 4
\setdots <1pt>
\ellipticalarc axes ratio 3:1.1  360 degrees from 2.5 0 center at 0 0
\ellipticalarc axes ratio 3:1.1  360 degrees from 2 0 center at 0 0
\ellipticalarc axes ratio 3:1.1  360 degrees from 1.5 0 center at 0 0
\ellipticalarc axes ratio 3:1.1  360 degrees from 1 0 center at 0 0
\ellipticalarc axes ratio 3:1.1  360 degrees from .5 0 center at 0 0
\setsolid
\ellipticalarc  axes ratio 3:1.1  360 degrees from 3 0 center at 0 0
\ellipticalarc axes ratio 3:1.1  360 degrees from 3 -4 center at 0 -4
\putrule from -3 -4 to -3 4
\putrule from 3 -4 to 3 4
\putrule from -1 -5 to -1 3
\putrule from 1 -3 to 1 5 
\put {$E_1$} at 2 -1.5
\put {$E$} at -1.6 2
\put {$Y$} at 3.6 3
\setlinear \plot -1 3  1 5 /
           \plot -1 -5  1 -3 /
  \plot -1 -1  1 1 /
  \plot -1 -.95  1 1.05 /
  \plot -1 -1.05  1 .95 /
\setlinear
\setshadegrid span <2truept>
\hshade -5 -1 -1   -3 -1 1    3 -1 1   5 1 1  /   
\endpicture}
\vskip 6truemm

As usual we denote $K_Y = h^\ast K_X + (\nu_1 - 1)E_1$ and $h^\ast (ND) = NE +
N_1 E_1$ for $N \in \Bbb Q$, $N > 0$.  One can verify that $\nu_1 = d + 1 - a$
and $N_1 = N$.

To compute $z_X (ND,s)$ we need the Euler characteristics of the varieties
$Y^\circ, E^\circ, E^\circ_1$ and $E \cap E_1$ (which stratify $Y$).
Since $X$ and $D$ are affine cones we have that
$$\chi(E^\circ) = \chi (D \setminus \{ 0 \}) = 0 \quad \text { and } \quad
\chi(Y^\circ) = \chi(X \setminus D)  = 0.$$
Now $E_1$ is a nonsingular hypersurface of degree $a$ in $\Bbb P^d$, yielding
$$\chi(E_1) = (1-a)\big(\frac{(1-a)^d-1}{a}\big) + d$$
(see for example [Hirz]).  And because $D$ was chosen to be general we have
moreover that $E \cap E_1$ is a nonsingular hypersurface of degree $a$ in $\Bbb
P^{d-1}$; so
$$\chi(E \cap E_1) = (1-a) \big(\frac{(1-a)^{d-1}-1}{a}\big) + d-1.$$
Then finally $\chi(E^\circ_1) = \chi (E_1) - \chi(E \cap E_1) = -(1-a)^d + 1$
and
$$\align
z_X(ND,s) = & z_0(ND,s) = \frac{\chi(E^\circ_1)}{\nu_1 + sN_1} + \frac{\chi(E
\cap E_1)}{(\nu_1 + sN_1)(1 + sN)} \\
= & \frac{-(1-a)^d + 1}{d+1-a+sN} +
\frac{(1-a)\big(\frac{(1-a)^{d-1}-1}{a}\big)+d-1}{(d+1-a+sN)(1+sN)} \\
= & \frac{(1-a)\big(\frac{(1-a)^d-1}{a}\big) + d + s\big(1 - (1-a)^d\big)N}{(d+1-a+sN)(1+sN)} .
\endalign
$$
A (not very exciting) calculation shows that there is no cancellation in this
expression, except when $d=3$ and $a = 2$ or $3$, in which case $z_X(ND,s)$ is
$$\frac{2}{1+sN} \quad \text{and } \quad \frac{9}{1+sN},$$
respectively.  Taking limits we obtain
$$e_X(ND) = \frac{(1-a)\big(\frac{(1-a)^d-1}{a}\big)+d+\big(1-(1-a)^d\big)N}{(d+1-a+N)(1+N)} \qquad
\text{ if }  d+1+N \ne a$$
and 
$$e_X\big((X,ND)\big) =
\frac{(1-a)\big(\frac{(1-a)^d-1}{a}\big)+d+\big((1-a)^d-1\big)N}{(d+1-a-N)(1-N)} \qquad \text{
if } 
\cases  d + 1 \ne a + N \\ N \ne 1. \endcases$$ 

\bigskip
\bigskip

\head
5. Zeta functions associated to divisors and differential forms
\endhead
\bigskip
\noindent
{\bf 5.1.}  In the $p$--adic theory of Igusa's local zeta functions one also
associates this invariant to both polynomials and differential forms, see e.g. [L, III3.5]. Let $f
\in \Bbb Q_p[x] = \Bbb Q_p[x_1,\cdots,x_d]$ and $w \in \Omega^d_{\Bbb A^d}$,
i.e. $w = g dx$ where $g \in \Bbb Q_p[x]$ and $dx = dx_1 \wedge \cdots \wedge
dx_d$.  Then, with the notation of 1.3.1,
$$Z_p(f,w,s) := \int_{\Bbb Z^d_p} |f(x)|^s |g(x)||dx|.$$
With the notation of 1.4 let $\nu^\prime_i - 1$ be the multiplicity of $E_i$ in
the divisor of $h^\ast w$.  Then (for $f \in \Bbb Q[x]$) the same formula as in 1.4.1 is valid when we
replace $\nu_i$ by $\nu^\prime_i$.

We also want to introduce this notion on the level of the Grothendieck ring of
algebraic varieties as in 1.3. Our motivation in this paper is that we will use
it to construct on a $\Bbb Q$--Gorenstein variety $X$ an invariant ${\Cal
Z}_W(X,D,s)$,
generalizing $Z_W(X,D,s)$ in 3.3, on the level of the Grothendieck ring.
Furthermore we will need this notion in future work.
\bigskip
\noindent
{\bf 5.2.} We fix notations for this section.  Let $X$ be an irreducible
nonsingular variety of dimension $d$ and $W$ a subvariety of $X$.  Let $D$ be
an effective divisor on $X$ and $J \subset \Omega^d_X$ an invertible subsheaf
of the sheaf of regular differential $d$--forms $\Omega^d_X$ on $X$.

We will only consider the situation where supp $J \subset
\supp D$; we motivate this below.
\bigskip
\noindent
{\bf 5.3.} First we rephrase the definition of $\Cal Z_W(D,s)$ in terms of the
{\it
motivic volume} $\mu$ of [DL3, 3.2] or [DL4].   Denote by $\Cal C$ the family
of subsets of $\Cal L(X)$ of the form $\pi^{-1}_n A_n$ for some $n \in \Bbb N$
and constructible subset $A_n$ of $\Cal L_n(X)$.  We call these {\it
cylindrical} subsets as in [B2] or [DL4].  There exists a unique additive
measure $\mu : \Cal C \rightarrow \Cal M_L$ satisfying $\mu(\pi^{-1}_n A_n) =
\frac{[A_n]}{L^{(n+1)d}}$ for $A_n$ as above.  (In fact this map is denoted by
$\tilde \mu$ in [DL3] and there
$\mu$ is a map from the more complicated family of semi--algebraic subsets of
$\Cal L(X)$ to $\hat M$.) 
For $A$ in $\Cal C$ and $\alpha : A \rightarrow \Bbb N$ a bounded function with
cylindrical fibres one defines the integral
$$\int_A L^{-\alpha} d \mu := \sum_{n \in \Bbb N} L^{-n} \mu(\alpha^{-1} \{n
\}) \in \Cal M_L. \tag 5.3.1$$
Now re--examining the definition of $\Cal Z_W(D,s)$ in 1.3 we have, with the
notation introduced there, that $\mu(Y_{n,D,W}) = [X_{n,D,W}] L^{-(n+1)d}$ and
hence
$$\Cal Z_W(D,s) = \sum_{n \in \Bbb N} \mu (Y_{n,D,W})L^{-ns} \in \Cal
M_L[[L^{-s}]].$$
\bigskip
\noindent
{\bf 5.4.} The following construction is a special case of [DL3, 3.5].  To the
sheaf $J$ is associated as follows a measure $\mu_J$ on $\Cal C$, such that
$\mu_{\Omega^d_X} = \mu$.

For $P \in X$ let $dx$ and $g_P dx$ be local generators of $\Omega^d_X$ and $J$,
respectively, around $P$.  Denote then by $\ord_t J : \Cal L(X) \rightarrow
\Bbb N \cup \{ \infty \}$ the function assigning to $\varphi$ in $\Cal L(X)$
the order of the power series given by $g_{\pi_0(\varphi)} \circ \varphi$.  For $A$ in $\Cal C$
we define 
$$\mu_J(A) := \int_A L^{-\ord_t J} d \mu = \sum_{\ell \in \Bbb N} L^{-\ell} \mu
(A \cap \{ \ord_t J = \ell \}).$$
Indeed the sets $\{ \ord_t J = \ell \}$ are cylindrical.  For arbitary $A$ the
right hand side above is only defined as an element in $\hat \Cal M$; however
we will only consider sets $A$ for which the sum over $\ell$ is finite and then
$\mu_J(A) \in \Cal M_L$.  Replacing $\mu$ by $\mu_J$ we can consider analogous
integrals as in (5.3.1).

The following change of variables formula is a special case of [DL3, 3.5.2].
(It follows immediately from [DL3, 3.3] of which Theorem 1.10 is a special
case.)
\bigskip
\proclaim
{5.4.1. Proposition}  Let $X^\prime$ be another irreducible smooth variety
and $\rho : X^\prime \rightarrow X$ a proper birational morphism.  For $A$ in
$\Cal C$ and $\alpha : A \rightarrow \Bbb N$ a bounded function with
cylindrical fibres we have that
$$\int_A L^{-\alpha} d \mu_J = \int_{\rho^{-1}A} L^{-\alpha \circ \rho} d
\mu_{\rho^\ast J}.$$
\endproclaim
\bigskip
\proclaim
{5.5. Definition} \rm To the data of 5.2 we associate the {\it motivic zeta
function}
$$\aligned
\Cal Z_W(D,J,s) & = \Cal Z_W(X,D,J,s) := \sum_{n \in \Bbb N}
\mu_J(Y_{n,D,W})L^{-ns} \\
& = \sum_{n \in \Bbb N} \big(\sum_{\ell \in \Bbb N} L^{-\ell} \mu (Y_{n,D,W} \cap
\{ \ord_t J = \ell \})\big)L^{-ns} \in \Cal M_L [[L^{-s}]].
\endaligned
$$
\endproclaim
\noindent
We explain why the sum over $\ell$ is finite.  For any fixed $n$ we have that
$Y_{n,D,W} = \coprod_{\ell \in \Bbb N \cup \{ \infty \}} (Y_{n,D,W} \cap \{
\ord_t J = \ell \})$.
But our condition $\supp J \subset \supp D$ implies that
$$\{ \ord_t J = \infty \} = \Cal L (\supp J) \subset \Cal L(\supp D) = \{ \ord_t
D = \infty \}\, ,$$
hence we have that $Y_{n,D,W} \cap \{ \ord_t J  = \infty \} = \emptyset$ and
so
$Y_{n,D,W}$ is the countable union of the {\it cylindrical} sets $Y_{n,D,W} \cap
\{ \ord_t J
= \ell \}, \ell \in \Bbb N$.  Then this union is finite by [B2, Theorem 6.6]. 
\bigskip
\proclaim
{5.6. Theorem}  Let $X^\prime$ be another irreducible smooth variety and $\rho
: X^\prime \rightarrow X$ a proper birational morphism.  Then
$$\Cal Z_W (X,D,J,s)  = \Cal Z_{\rho^{-1}W} (X^\prime, \rho^\ast D, \rho^\ast
J, s).$$
\endproclaim
\bigskip
\demo {Proof}  This is a consequence of Proposition 5.4.1. \qed
\enddemo
\bigskip
\proclaim
{5.7. Theorem}  Let $h : Y \rightarrow X$ be a log resolution of $\supp D$.
  Denote as usual the irreducible components of $h^{-1}(\supp D)$ by
$E_i, i \in T$.  We set $h^\ast D = \sum_{i \in T} N_i E_i$ and $\div(h^\ast w)
= \sum_{i \in T} (\nu^\prime_i - 1)E_i$, where $w$ is a local generator of $J$.
Then
$$\Cal Z_W(D,J,s) = L^{-d} \sum_{I \subset T} [E^\circ_I \cap h^{-1}W] \prod_{i \in I}
\frac{L-1}{L^{\nu^\prime_i + sN_i} - 1}.$$
\endproclaim
\bigskip
\noindent
{\it Remark.}  Let as in 1.4 $dx$ be a local generator of $\Omega^d_X$ and
$\div(h^\ast dx) = \sum_{i \in T} (\nu_i - 1)E_i$.  Say $w = g dx$ and
$\div(h^\ast g) = \sum_{i \in T} M_i E_i$.  Then $\nu^\prime_i = \nu_i + M_i$
for $i \in T$.
\bigskip
\noindent
\demo{Proof}  One can adapt the proof of [DL2, Theorem 2.2.1] completely to
this more general setting with the sheaf $J$. \qed
\enddemo
\bigskip
\noindent
{\bf 5.8.} The notion introduced above is sufficient to introduce zeta
functions on the level of the Grothendieck ring for Gorenstein varieties.  To
cover the  case of $\Bbb Q$--Gorenstein varieties we need `sheaves of
multivalued differential forms'.  We briefly describe this generalization. 

Now let $J \subset (\Omega^d_X)^{\otimes m}$ be an invertible subsheaf of the
$m$--fold tensor product of $\Omega^d_X$, still satisfying $\supp J \subset
\supp D$.  We define
$$\mu_{J^{1/m}} (A) := \int_A L^{-\frac{\ord_tJ}{m}} d \mu = \sum_{\ell \in
\Bbb N} L^{-\ell/m} \mu (A \cap \{ \ord_t J = \ell \}) \in \Cal M_L [L^{1/m}]$$
for the sets $A$ in $\Cal C$ for which the last sum is finite.  Then the
motivic zeta function is
$$\aligned 
\Cal Z_W(D,J^{1/m},s) & = \Cal Z_W (X,D,J^{1/m},s) \\
& := \sum_{n \in \Bbb N} \mu_{J^{1/m}} (Y_{n,D,W})L^{-ns} \in \Cal
M_L[L^{1/m}][[L^{-s}]].
\endaligned
$$
Theorem 5.7 easily generalizes to this setting, but now the $\nu_i^\prime \in
\frac 1m
(\Bbb N \setminus \{ 0 \})$.
\bigskip
\noindent
{\bf 5.9.} Finally as in 1.5 we can generalize further to $\Bbb Q$--divisors.
Now if
$D$ is an effective $\Bbb Q$--divisor on $X$, such that $rD$ is a divisor for an
$r \in \Bbb N \setminus \{ 0 \}$, we define $\Cal Z_W(D,J^{1/m},s) := \Cal Z_W
(rD,J^{1/m},s/r)$.  Again Theorem 5.7 generalizes, with now the $N_i \in \frac
1r (\Bbb N \setminus \{ 0 \})$.
\bigskip
\bigskip
\head
6. Singular varieties; on the level of the Grothendieck ring
\endhead
\bigskip
\noindent
{\bf 6.1.}  In this section we generalize the zeta function of 3.3 to the level
of the Grothendieck ring.  In order to focus on the main idea we first treat
the essential case, being an effective Cartier divisor $D$ on a Gorenstein
variety $X$. For a normal variety $V$ we denote its canonical sheaf (corresponding to $K_V$) by $\omega_V$; we have that $\omega_V$ is invertible or $\omega_V^{\otimes m}$ is invertible for some $m \in \Bbb N \setminus \{0\}$ precisely when $V$ is Gorenstein or $\Bbb Q$--Gorenstein, respectively.

Also in the sequel $\Cal I(F)$ denotes the sheaf of ideals associated
to an effective divisor $F$ on a nonsingular variety.
\bigskip
\proclaim
{6.2. Definition -- Proposition}  Let $X$ be a Gorenstein variety of
dimension $d$ and $W$ a subvariety of $X$.  Let $D$ be an effective Cartier
divisor on $X$ such that $X_{\sing} \subset \supp D$.

(i) The {\it motivic zeta function}
$$\Cal Z_W(D,s) = \Cal Z_W(X,D,s) := \Cal Z_{h^{-1}W} (Y,h^\ast D, h^\ast
\omega_X \otimes \Cal
I(ah^\ast D), s-a)$$
where $h : Y \rightarrow X$ is a log resolution of $\supp D$ and $a \in \Bbb N,
a >> 0$.

(ii) Let $h : Y \rightarrow X$ be a log resolution of $\supp D$.  With the
notation of 3.2 we have that
$$\Cal Z_W(D,s) = L^{-d} \sum_{I \subset T} [E^\circ_I \cap h^{-1} W] \prod_{i
\in I} \frac{L-1}{L^{\nu_i + sN_i} - 1}.$$
\endproclaim

\demo
{Proof} We first explain the right hand side of our definition.  Since $X_{\sing} \subset \supp D$ we have that
$\supp(h^\ast \omega_X) \subset \supp(h^\ast D)$, yielding for 
$a >> 0$ that $h^\ast \omega_X \otimes \Cal I(ah^\ast D)$ is an invertible
subsheaf of $\Omega^d_Y$.  So to this sheaf and the
effective divisor $h^\ast D$ we can associate the motivic zeta function of 5.5.
The substitution `$s-a$ instead of $s$'  means replacing the variable
$L^{-s}$
by
$L^a(L^{-s})$. 

Now we show independency of the chosen resolution; it is sufficient to consider
another log resolution $h^\prime : Y^\prime \rightarrow X$ that factors as
$h^\prime : Y^\prime \overset \varphi \to \rightarrow Y \overset h \to
\rightarrow X$.  By Theorem 5.6 we indeed have that
$$\aligned
 \Cal Z_{h^{-1}W} & (Y, h^\ast D, h^\ast \omega_X \otimes \Cal I(ah^\ast
 D),s-a) \\
& =  
\Cal Z_{\varphi^{-1}(h^{-1}W)} (Y^\prime, \varphi^\ast h^\ast D, \varphi^\ast
(h^\ast \omega_X) \otimes \varphi^\ast(\Cal I(ah^\ast D)),s-a) \\
& = \Cal Z_{h^{\prime -1}W} (Y^\prime, h^{\prime \ast} D, h^{\prime \ast}
\omega_X \otimes \Cal I (ah^{\prime \ast} D), s-a).
\endaligned
$$
Let $\eta$ and $f$ be local generators of $\omega_X$ and $\Cal I(D)$,
respectively.  Then $(h^\ast f)^a (h^\ast \eta)$ is a local generator of
$h^\ast \omega_X \otimes \Cal I(ah^\ast D)$ and its divisor of zeroes is
$\sum_{i \in T} ((\nu_i - 1)  + aN_i)E_i$.  Hence Theorem 5.7 (with $h=Id_Y$) yields the stated
formula for $\Cal Z_W(D,s)$, which also proves independency of the number $a$.
\qed
\enddemo
\bigskip
\noindent
{\bf  6.3.}  Now let $X$ be $\Bbb Q$--Gorenstein and say that $mK_X$ is Cartier
for some $m \in \Bbb N \setminus \{ 0 \}$.   We define $\Cal Z_W (X,D,s)$ just
as in 6.2, interpreting the expression $h^\ast \omega_X \otimes \Cal I(ah^\ast
D)$ as an abbreviation of $(h^\ast (\omega_X^{\otimes m}) \otimes \Cal
I(mah^\ast D))^{1/m}$ (see 5.8).  Now $\Cal Z_W(X,D,s)$ lives in a localization
of $\Cal M_L[L^{1/m}][L^{-s}]$ and is given by the same formula as in 6.2 (with
now the $\nu_i \in \Bbb Q$).

When $D$ is an effective $\Bbb Q$--Cartier divisor we set as usual $\Cal
Z_W(D,s) := \Cal Z_W (rD,s/r)$ if $rD$ is Cartier for an $r \in \Bbb N
\setminus
\{ 0 \}$.  Then in full generality we have the following.
\bigskip
\proclaim
{6.4. Definition -- Proposition}  Let $X$ be a $\Bbb Q$--Gorenstein variety of
dimension $d$ and $W$ a subvariety of $X$.  Let $D$ be an effective $\Bbb
Q$--Cartier divisor on $X$ (with $rD$ Cartier for an $r \in \Bbb N \setminus \{
0
\})$ such that $X_{\sing} \subset \supp D$.  The {\it motivic zeta function}
$$
\Cal Z_W(D,s) = \Cal Z_W(X,D,s) := \Cal Z_{h^{-1}W} (Y,h^\ast(rD),h^\ast
\omega_X
\otimes \Cal I(ar h^\ast D), s/r - a)$$ 
where $h : Y \rightarrow X$ is a log resolution of $\supp D$ and $a \in \Bbb
N, a >> 0$.  We have the same formula as in 6.2.
\endproclaim
\bigskip
\noindent
Of course $\Cal Z_W(D,s)$ specializes to the zeta function $Z_W(D,s)$ of
3.3.
\bigskip
\noindent
{\bf 6.5.} Finally we consider for arbitrary $\Bbb Q$--Gorenstein varieties $X$
`Kontsevich' invariants $\Cal E_W(D)$ and $\Cal E_W \bigl((X,D)\bigr)$ on the level of the
Grothendieck ring, which specialize to $E_W(D)$ and $E_W \bigl((X,D)\bigr)$ of 3.6 and
3.10, respectively.
Notice first that in 6.4 we have, by the formula for $\Cal Z_W (D,s)$ in terms of a log resolution, that it already belongs to the localization of a polynomial ring $\Cal M_L[L^{1/r}] [L^{-s/r}]$ with respect to $(1-L^{-\alpha - \beta s})_{\alpha \in \Bbb Q, \beta \in \Bbb Q_{>0}}$.
Morally we again take limits for $s \rightarrow 1$ and $s \rightarrow -1$ to define  $\Cal E_W(D)$ and $\Cal E_W \bigl((X,D)\bigr)$, respectively.
\bigskip
\noindent
\proclaim
{6.6. Definition} \rm  Let $X$ be a $\Bbb Q$--Gorenstein variety of
dimension $d$ and $W$ a subvariety of $X$.  Let $D$ be an effective $\Bbb
Q$--Cartier divisor on $X$ such that $X_{\sing} \subset \supp D$.  Take $r
\in \Bbb N \setminus \{ 0 \}$ with $rK_X$ and $rD$ Cartier.

(i) If $\Cal Z_W (D,s)$ belongs to the localization of $\Cal M_L[L^{1/r}] [L^{-s/r}]$ with respect to    \linebreak
$(1-L^{-\alpha - \beta s})_{\alpha \in \Bbb Q, \beta \in \Bbb Q_{>0}, \alpha+\beta\neq 0}$, then we put
$$
\Cal E_W(D) = \Cal E_W(X,D) := \Cal Z_W(D,s)|_{s=1} \, .
$$
\noindent
Otherwise we put $\Cal E_W(D)= \Cal E_W(X,D) := \infty$.

(ii) If $\Cal Z_W (D,s)$ belongs to the localization of $\Cal M_L[L^{1/r}] [L^{-s/r}]$ with respect to   \linebreak
$(1-L^{-\alpha - \beta s})_{\alpha \in \Bbb Q, \beta \in \Bbb Q_{>0}, \alpha\neq\beta}$, then we put
$$
\Cal E_W \bigl((X,D)\bigr) := \Cal Z_W(D,s)|_{s=-1} \, .
$$
Otherwise we put $\Cal E_W \bigl((X,D)\bigr) := \infty$.

Here the evaluations $s=1$ and $s=-1$ mean substituting the variable
$L^{-s/r}$ by $L^{-1/r}$ and $L^{1/r}$, respectively, yielding a well defined
element in $\Cal{\hat M}[L^{1/r}]$.
\endproclaim
\bigskip
\proclaim
{6.7. Proposition}  Consider the same data as is 6.6.

(i) Suppose there is a log resolution $h : Y \rightarrow  X$ of $\supp D$ for
which $\nu_i + N_i \ne 0$ for all $i \in T$ (using the notation of 3.2).  Then
$$
\Cal E_W(D) = L^{-d} \sum_{I \subset T} [E^\circ_I \cap h^{-1} W] \prod_{i \in I}
\frac{L-1}{L^{\nu_i+N_i}-1} .
$$

(ii) Suppose there is a log resolution $h : Y  \rightarrow X$ of $\supp D$ for
which all log discrepancies $a_i, i \in T$, with respect to the pair $(X,D)$
satisfy $a_i \ne 0$ (using the notation of 3.2).  Then
$$
\Cal E_W \bigl((X,D)\bigr) = L^{-d} \sum_{I \subset T} [E^\circ_I \cap h^{-1} W] \prod_{i \in I}
\frac{L-1}{L^{a_i}-1} .
$$
\endproclaim
\bigskip
\bigskip
\head
Appendix
\endhead
\bigskip
\noindent
{\bf A1.} Let in this appendix $X$ be a smooth irreducible variety of dimension
$d$ and $W$ a subvariety of $X$.

In 1.7 -- 1.10 we described the Kontsevich invariant $\Cal E_W(D) \in \hat {\Cal
M}$, associated to an effective divisor $D$ on $X$, and we mentioned its
important properties.  Here we will generalize this notion to effective $\Bbb
Q$--divisors; if $rD$ is a divisor for an $r \in \Bbb N \setminus \{ 0 \}$ we
obtain an invariant $\Cal E_W(D)$ in a finite extension $\hat \Cal M[L^{1/r}]$ of
$\hat \Cal M$,  and we treat analogous properties.  We also introduce this
invariant for a $\Bbb Q$--divisor $D = \sum_i d_i D_i$ (with the $D_i$ irreducible)
such that all $d_i > -1$ and $\supp D = \cup_i D_i$ is a divisor with strict
normal crossings.  This is used in 3.14.
\bigskip
\noindent
{\bf A2.} First we describe the ring $\hat \Cal M[L^{1/r}]$.  Consider the
integral ring extension $\Cal M_L \hookrightarrow \Cal M_L [L^{1/r}] := {\Cal
M_L[X]
\over (X^r-L)}$, where $L^{1/r}$ is the class of $X$ in this quotient.
Each element $a \in \Cal M_L[L^{1/r}]$ has a unique expression of the form $a =
\sum^{r-1}_{i=0} a_i L^{i/r}$ or $a = \sum^{r-1}_{i=0} a^\prime_i L^{-i/r}$
with $a_i, a^\prime_i \in \Cal M_L$.

We extend the decreasing filtration $(F_m)_{m \in \Bbb Z}$ on $\Cal M_L$,
introduced in 1.7, to the ring $\Cal M_L [L^{1/r}]$. Let $F^\prime_m, m \in
\Bbb Z$, be the subgroup of $\Cal M_L[L^{1/r}]$ generated by
$$\big\{ \sum^{r-1}_{i=0} {[A_i] \over L^{n_i}} L^{-i/r} | \dim A_i - n_i \leq -m
\text{ for } i = 0, \dots , r-1 \big\}.$$ 
(So indeed $F_m = \Cal M_L \cap F^\prime_m$.)  We take the completion $\hat
\Cal M^\prime$ of $\Cal M_L[L^{1/r}]$ with respect to this filtration
$(F^\prime_m)_{m \in \Bbb Z}$; then we have an injection $\hat \Cal M
\hookrightarrow \hat \Cal M^\prime$.

One can verify that $\hat \Cal M' \cong \hat \Cal M[L^{1/r}]$, where the right
hand side can be interpreted either as the subring of $\hat \Cal M^\prime$ generated by
$\hat \Cal M$ and $L^{1/r}$, or as ${\hat \Cal M[X] \over (X^r - L)}$.
\bigskip
\noindent
{\bf A3.}  We will use the following notation.  Let $D$ be a prime divisor on
$X$.  Then $\ord_t D : \Cal L(X) \rightarrow \Bbb N \cup \{ \infty \}$ assigns
to $\varphi \in \Cal L(X)$ the order of the power series in $t$ given by $f \circ \varphi$, where $f$
is a local equation of $D$ at $\pi_0(\varphi)$.  For a $\Bbb Q$--divisor $D  =
\sum_i d_i D_i$ (with the $D_i$ prime divisors) we then define $\ord_t D : \Cal
L(X) \rightarrow \Bbb Q \cup \{ \infty \}$ by $\ord_t D := \sum_i d_i \ord_t
D_i$.  
\bigskip
\noindent
{\bf A4. Definition.}  Let $D$ be a $\Bbb Q$--divisor on $X$ and $r \in \Bbb N
\setminus \{ 0 \}$ such that $rD$ is a divisor.  

($i$) If $D$ is effective we
define for $n \in \Bbb N$ the subscheme $Y_{n,D,w}$ of $\Cal L(X)$ and the
subscheme $X_{n,D,w}$ of $\Cal L_n(X)$ as in 1.3 with only the following
adaptation : now $f$ is a local equation of the divisor $rD$ (instead of $D$).
Then we set
$$\Cal E_W(D) = \Cal E_W(X,D) := \sum_{n \in \Bbb N}
\frac{[X_{n,D,W}]}{L^{(n+1)d}}L^{-n/r} \in \hat \Cal M[L^{1/r}].$$

\noindent
In terms of the motivic volume $\mu$ of 5.3 we can describe $\Cal E_W(D)$ as
$$\Cal E_W(D) = \int_{\pi^{-1}_0 W} L^{-\ord_t D} d \mu := \sum_{n \in \Bbb N}
\mu (\pi^{-1}_0 W \cap \{ \ord_t D =  \frac nr \})L^{-n/r}.$$

($ii$) In general we say that $\ord_t D : \Cal L(X) \rightarrow \frac 1r \Bbb Z
\cup \{ \infty \}$ is integrable on $\pi^{-1}_0 W$ if
$$\int_{\pi^{-1}_0 W} L^{-\ord_t D} d \mu := \sum_{n \in \Bbb Z} \mu(\pi^{-1}_0
W \cap \{ \ord_t D = \frac nr \})L^{-n/r}$$ 
converges in $\hat \Cal M[L^{1/r}]$; we then denote this invariant again
by $\Cal E_W(D)$.
\bigskip
\noindent
{\bf A5.}   An important case of this last definition occurs when $D =
\sum^k_{i=1} d_i D_i$ with the $D_i$ irreducible, all $d_i > -1$, and supp
$D = \cup^k_{i=1} D_i$ a divisor with strict normal crossings.  For $J \subset
\{ 1, \cdots , k \}$ denote $D^\circ_J := (\cap_{j \in J} D_j)
\setminus (\cup_{\ell \notin J} D_\ell)$ and $M_J := \{ (m_1,\cdots,m_k) \in \Bbb N^k \mid$ $m_j > 0 \Leftrightarrow j \in J \}$.  Then one can compute that   
$$\int_{\pi^{-1}_0 W} L^{-\ord_t D} d \mu = L^{-d} \sum_{J \subset \{ 1,
\cdots , k \} } (L-1)^{|J|} [D_J^\circ \cap W] \sum_{(m_1,\cdots,m_k) \in M_J}
L^{-\sum_{j \in J} (d_j+1)m_j} ,$$
which converges in $\hat \Cal M[L^{1/r}]$ since all $d_j + 1 > 0$.  See
[B2, Theorem 6.28] and [C, Theorem 1.17].
\bigskip
\proclaim
{\bf A6. Theorem}  Let also $X^\prime$ be a smooth irreducible variety and
$\rho : X^\prime \rightarrow X$ a proper birational morphism.  Let $D$ be a
$\Bbb Q$--divisor on $X$.  Then $\ord_t D$ is integrable on $\pi^{-1}_0 W$ if
and only if $\ord_t(\rho^\ast D + K_{X^\prime|X})$ is integrable on $\pi^{-1}_0
(\rho^{-1} W)$; and in this case
$$\Cal E_W(X,D) = \Cal E_{\rho^{-1}W} (X^\prime, \rho^\ast D + K_{X^\prime|X}).$$
\endproclaim
\bigskip
\demo
{Proof}  The proof of [DL3, Lemma 3.3], based on the crucial and difficult
[DL3, Lemma 3.4], can be adapted to this setting.  See [B2, Theorem 6.27] for an
analogous statement and proof when $W = X$.  We also remark that when $D$ is an
{\it effective} $\Bbb Q$--divisor (implying that both functions are integrable), then
one can prove the stated equality {\it using} the equality in [DL3, Lemma 3.3]. \qed
\enddemo
\bigskip
\proclaim
{A7. Theorem} Let $D$ be a $\Bbb Q$--divisor on $X$ (with $rD$ a divisor for an
$r \in \Bbb N \setminus \{ 0 \}$) such that $\ord_t D$ is integrable on
$\pi^{-1}_0 W$.  Using the notation of 1.4 we have the following formula for
$\Cal E_W(D)$ in terms of a log resolution $h : Y \rightarrow X$ of supp $D$ :
$$\Cal E_W(D) = L^{-d} \sum_{I \subset T} [E^\circ_I \cap h^{-1} W] \prod_{i
\in I} \frac{L-1}{L^{\nu_i + N_i}-1} \quad \text{ in } \quad \hat \Cal
M[L^{1/r}].$$
In particular $\Cal E_W(D)$ belongs to the image of $\Cal M_L[(1 -
L^{-n/r})^{-1}]_{n \in \Bbb N \setminus \{ 0 \}}$ in $\hat \Cal M[L^{1/r}]$.
\endproclaim
\bigskip
\demo
{Proof}  This follows from A5 and Theorem A6.  One can also adapt [DL3, (6.5)].
\qed
\enddemo


\bigskip
\bigskip

\Refs

\ref \key AKMW
\by D. Abramovich, K. Karu, K. Matsuki and J. W{\l}odarczyk
\paper Torification and factorization of birational maps
\jour math.AG/9904135
\vol 
\yr 1999
\pages 
\endref

\ref \key B1
\by V\. Batyrev
\paper Birational Calabi--Yau $n$--folds have equal Betti numbers
\inbook New Trends in Algebraic geometry, Euroconference on Algebraic Geometry (Warwick 1996), London Math. Soc. Lecture Note Ser. 264, K. Hulek et al Ed. 
\vol 
\publ CUP
\yr 1999
\pages 1--11
\endref

\ref \key B2
\by V\. Batyrev
\paper Stringy Hodge numbers of varieties with Gorenstein canonical singularities 
\jour Proc. Taniguchi Symposium 1997, In \lq Integrable Systems and Algebraic Geometry, Kobe/ \linebreak Kyoto 1997\rq, World Sci. Publ.
\vol 
\yr 1999
\pages 1--32 
\endref

\ref \key B3
\by V\. Batyrev
\paper Non--Archimedian integrals and stringy Euler numbers of log terminal pairs
\jour J. Europ. Math. Soc.
\vol 1
\yr 1999
\pages 5--33
\endref

\ref \key BLR
\by S\. Bosch, W\. L\"utkebohmert and M\. Raynaud
\book N\'eron Models
\bookinfo Ergeb. Math. Grenzgeb. (3) 21
\publ Springer Verlag, Berlin
\yr 1990
\endref

\ref \key C
\by A\. Craw
\paper An introduction to motivic integration
\jour math.AG/9911179
\vol   
\yr 1999
\pages 
\endref

\ref \key D1
\by J\. Denef
\paper On the degree of Igusa's local zeta function
\jour Amer. J. Math.
\vol 109
\yr 1987
\pages 991--1008
\endref

\ref \key D2
\by J\. Denef
\paper Report on Igusa's local zeta function
\jour Ast\'erisque
\paperinfo S\'em. Bourbaki 741
\vol 201/203
\yr 1991
\pages 359--386
\endref

\ref \key DL1
\by J\. Denef and F\. Loeser
\paper Caract\'eristiques d'Euler--Poincar\'e, fonctions zeta locales, et
modifications analytiques
\jour J. Amer. Math. Soc.
\vol 5
\yr 1992
\pages 705--720
\endref

\ref \key DL2
\by J\. Denef and F\. Loeser
\paper Motivic Igusa zeta functions
\jour J. Alg. Geom.
\vol 7
\yr 1998
\pages 505--537
\endref

\ref \key DL3
\by J\. Denef and F\. Loeser
\paper Germs of arcs on singular algebraic varieties and motivic integration
\jour Invent. Math.
\vol 135
\yr 1999
\pages 201--232
\endref

\ref \key DL4
\by J\. Denef and F\. Loeser
\paper Motivic integration, quotient singularities and the McKay correspondence
\jour preprint
\vol 
\yr 1999
\pages 
\endref

\ref \key Hiro
\by H\. Hironaka
\paper Resolution of singularities of an algebraic variety over a field of
       characteristic zero
\jour Ann. Math.
\vol 79
\yr 1964
\pages 109--326
\endref

\ref\key Hirz
\by F\. Hirzebruch
\book Topological methods in algebraic geometry
\publ Springer Verlag
\yr 1966
\endref

\ref \key I
\by J\. Igusa
\paper Complex powers and asymptotic expansions I
\jour J. Reine Angew. Math.
\vol 268/269
\yr 1974
\pages 110--130
\moreref
\paper II
\jour ibid.
\vol 278/279
\yr 1975
\pages 307--321
\endref

\ref \key KM
\by J\. Koll\'ar and S\. Mori
\book Birational geometry of algebraic varieties
\bookinfo Cambridge Tracts in Mathematics 134
\publ Cambridge Univ. Press
\yr 1998
\endref

\ref \key KMM
\by Y\. Kawamata, K\. Matsuda and K\. Matsuki
\paper Introduction to the Minimal Model Program
\jour Adv. Stud. Pure Math.
\paperinfo Algebraic Geometry, Sendai, T\. Oda ed., Kinokuniya
\vol 10
\yr 1987
\pages 283--360
\endref

\ref \key Kol
\by J\. Koll\'ar et al.
\paper Flips and abundance for algebraic threefolds
\jour Ast\'erisque
\paperinfo seminar Salt Lake City 1991
\vol 211
\yr 1992
\endref

\ref \key Kon
\by M\. Kontsevich
\paper 
\jour Lecture at Orsay (December 7, 1995)
\yr
\pages
\endref

\ref \key L
\by F\. Loeser
\paper Fonctions d'Igusa $p$--adiques et polyn\^omes de Bernstein
\jour Amer. J. Math.
\vol 110
\yr 1988
\pages 1--22
\endref

\ref \key M
\by D\. Mumford
\paper The topology of normal singularities of an algebraic surface and a
criterion for simplicity
\jour Publ. Math. I.H.E.S.
\vol 9
\yr 1961
\pages 5--22
\endref

\ref \key S
\by E. Szab\'o
\paper Divisorial log terminal singularities
\jour J. Math. Sci. Univ. Tokyo  
\vol 1
\yr 1994
\pages 631--639
\endref

\ref \key V1
\by W. Veys
\paper Determination of the poles of the topological zeta function for curves
\jour Manuscripta Math.
\vol 87
\yr 1995
\pages 435--448
\endref

\ref \key V2
\by W. Veys
\paper Zeta functions for curves and log canonical models
\jour Proc. London Math. Soc.
\vol 74  
\yr 1997
\pages 360--378
\endref


\ref \key V3
\by W. Veys
\paper The topological zeta function associated to a function on a normal
surface germ
\jour Topology
\vol 38
\yr 1999
\pages 439--456
\endref

\ref \key Wa
\by C.-L. Wang
\paper On the topology of birational minimal models
\jour J. Differential Geom. 
\vol 50
\yr 1998
\pages 129--146
\endref


\ref \key W{\l}
\by J. W{\l}odarczyk
\paper Combinatorial structures on toroidal varieties and a proof of the weak factorization theorem
\jour math.AG/9904076
\vol 
\yr 1999
\pages 
\endref


\endRefs
\enddocument